\documentclass[12pt,a4paper]{article}
\usepackage{textcomp}
\usepackage{listings}
\usepackage[space=true]{accsupp}
\newcommand{\copyablespace}{\BeginAccSupp{method=hex,unicode,ActualText=00A0}\EndAccSupp{}}
\usepackage[T1]{fontenc}
\usepackage[english]{babel}
\usepackage{hyperref}
\usepackage{graphicx}
\usepackage{amsmath}
\usepackage{amssymb}
\usepackage{mathtools}
\usepackage[table,dvipsnames]{xcolor}
\usepackage{fancyhdr}
\usepackage{longtable}
\usepackage{multicol}
\usepackage{enumerate}
\usepackage{enumitem}
\usepackage{multirow}
\usepackage{natbib}
\usepackage{caption}
\usepackage{subcaption}
\usepackage{xurl}
\usepackage[utf8]{inputenc}
\usepackage{tikz}
\usetikzlibrary{shapes}
\usepackage{pgfplots}
\usepackage{amsthm}

\pgfplotsset{compat=newest}

\usepackage{framed}
\usepackage{booktabs}
\usepackage{caption}
\usepackage{float}
\usepackage{titlesec}
\usepackage{capt-of}
\definecolor{darkgreen}{rgb}{0.1, 0.5, 0.2}
\usepackage{array}
\usepackage{arydshln}
\newtheorem{definition}{Definition}
\newtheorem{example}{Example}
\newtheorem{proposition}{Proposition}
\newtheorem{theorem}{Theorem}
\newtheorem{remark}{Remark}
\newtheorem{corollary}{Corollary}
\newtheorem{assumption}{Assumption}
\newtheorem{lemma}{Lemma}

\newcommand{\Longupdownarrow}{\Big\Updownarrow}

\usepackage{esvect}
\usetikzlibrary{intersections}
\usetikzlibrary{backgrounds}
\definecolor{forestgreen}{rgb}{0.13, 0.55, 0.13}
\usepackage{siunitx}

\title{Geometry of efficient weight vectors}
\author{Krist\'of \'Abele-Nagy$^a$, S\'andor Boz\'oki$^{b,a}$, Zsombor Sz\'adoczki$^{b,a,*}$}
\date{}
\begin{document}
\pagenumbering{arabic}

\maketitle
\begin{center}
$*$ Corresponding author, 1111 Kende u. 13-17., Budapest, Hungary;\\ Email: szadoczki.zsombor@sztaki.hun-ren.hu\\
\bigskip
$^{a}$ Department of Operations Research and Actuarial Sciences \\
Corvinus University of Budapest, 1093 Fővám tér 8., Budapest, Hungary \\
Email: kristof.abele-nagy@uni-corvinus.hu\\
\bigskip
$^{b}$ Research Group of Operations Research and Decision Systems, \\
Research Laboratory on Engineering \& Management Intelligence \\
HUN-REN Institute for Computer Science and Control (SZTAKI), 1111 Kende u. 13-17., Budapest, Hungary;\\ Email: szadoczki.zsombor@sztaki.hun-ren.hu, bozoki.sandor@sztaki.hun-ren.hu\\

\end{center}

\bigskip

\newpage

\begin{abstract}

\noindent
Pairwise comparison matrices and the weight vectors obtained from them are important concepts in multi-criteria decision making. A weight vector calculated from a pairwise comparison matrix is called Pareto efficient if the approximation of the matrix elements by the weight ratios cannot be improved for any element of the matrix without worsening it for another element. The aim of this paper is to show the geometrical properties of the set of Pareto efficient weight vectors for $4\times 4$ pairwise comparison matrices. We prove that the set of efficient weight vectors is a union of three tetrahedra, each determined by four weight vectors calculated from an incomplete submatrix of the pairwise comparison matrix that can be represented by a path spanning tree graph. It is shown that with suitable rearrangements the orientations of the $4$-cycles in the Blanquero-Carrizosa-Conde graphs for the efficient weight vectors are determined as well. The special cases (double perturbed, simple perturbed, consistent matrices) are discussed in the online appendices.

\end{abstract}

\noindent \textbf{Keywords}:  Pairwise comparison, Pareto optimality, Geometry, Multi-criteria decision making, Analytic Hierarchy Process, Efficient weight vector, Spanning tree

\bigskip 

\noindent \textbf{MSC2020}: 05C20, 15B48, 90B50, 91B06


\renewcommand{\baselinestretch}{1.24} \normalsize

\section{Introduction}
\label{sec:1}

Pairwise comparisons are popular in many fields, such as decision making, psychology, and sports \citep{Thurstone1927,Triantaphyllou2000,Csato2021}.
Pairwise comparison matrices, originally developed for the Analytic Hierarchy Process \citep{Saaty1977,Saaty1980}, are important tools of multi-criteria decision making. When faced with a decision problem that involves multiple criteria, it is often not feasible to give the importance weights of the criteria directly. Instead, comparing the importance of the criteria pairwise in a ratio form by recording the answers to the questions ‘How many times is Criterion $i$ more important than Criterion $j$?' is much more feasible. These answers can be arranged in a matrix, which is called a pairwise comparison matrix (PCM).
The next step is to obtain a weight vector from the pairwise comparison matrix. This weight vector represents the estimation of the real underlying preferences of the decision maker. Instead of criteria, this approach can be applied in the same way to the performance of alternatives with respect to a single criterion. In this paper, the latter narrative is followed.

A PCM provided by a decision maker usually contains some contradictions, a certain level of inconsistency can be observed in the data. This can be measured by several inconsistency indices proposed in the literature \citep{Brunelli2018}.

It is also worth mentioning that to estimate the decision maker's underlying preferences (calculate a weight vector), not all of the pairwise comparisons are necessary. When there are missing elements in a PCM, then it is called an incomplete pairwise comparison matrix \citep{Harker1987}. The structure of the known pairwise comparisons is often represented by undirected graphs \citep{Gass1998}.

Since the matrix elements are the judgments of the decision maker, the weight ratios (regardless of how the weight vector is obtained), which are the approximations of the judgments, should be as close as possible to the matrix elements. If no approximation can be improved by changing the weights for any matrix element without making the approximation worse in some other element (thus, there is no other weight vector that (Pareto) dominates this one), the weight vector is called (Pareto) efficient. This is a natural requirement of a weight vector, since if the approximation of judgments can be improved in a trivial way, that improvement should be made.

Unfortunately, even popular weight vector calculation methods can produce inefficient weight vectors \citep{Bozoki2014}. Thus, describing the set of efficient weight vectors is an important task. The efficiency of weight vectors has recently gained more attention resulting in several new findings and studies in the related literature. However, almost all of those are only using the algebraic approach, while this paper aims to contribute to the complementary geometric view, and connect the two approaches as well.

\cite{Blanquero} besides dealing with the efficiency of some of the most popular weight calculation techniques, connected the problem of efficient weight vectors to a directed graph, which was used widely later on. 

\cite{Abele2016} and \cite{Abele2018} examined the popular eigenvector weight calculation method, and found that it is efficient for simple perturbed and double perturbed PCMs, i.e., for those ones, which only differ from a consistent matrix by one or two element(s) (and their reciprocal values), respectively. This result was extended to a certain class of triple perturbed PCMs---which differ from the consistent case by three elements (and their reciprocals)---by \cite{Fernandes2022}. \cite{FurtadoJohnson2024b} showed different classes of PCMs for which the eigenvector method is (in)efficient.

The efficiency of so-called block perturbed PCMs was studied by \cite{FurtadoJohnson2024a}.  All efficient weight vectors of simple \cite{DaCruz2021} and double \cite{Furtado2023} perturbed PCMs were also described in an algebraic way in the literature.

\cite{Conde2010} used linear optimization problems to
derive efficient weight vectors and determined a consistency index with an associated weakly efficient weight vector, while \cite{BozokiFulop2018} also developed LP models and algorithms to improve an inefficient weight vector, and proved that the eigenvector method is weakly efficient.

The efficiency of weight vectors calculated from PCMs related to a real-world problem was examined by \cite{DulebaMoslem2019}.

\cite{FurtadoJohnson2024e} proved that the geometric mean of any collection of distinct columns of a PCM
is an efficient weight vector, which was then extended by \cite{FurtadoJohnson2024c}, who showed that any Hadamard weighted geometric mean of the columns of a PCM is efficient. \cite{FurtadoJohnson2024d} also examined the efficiency of the convex hull of the columns of certain PCMs.

\cite{FurtadoJohnson2024f} deal with PCMs for which all efficient weight vectors provide the same ranking of alternatives, while \cite{FurtadoJohnson-completeset-2025} provided an algebraic algorithm that determines inductively the complete set of efficient weight vectors for a PCM, as well as proved that the set of efficient weight vectors is piecewise linear connected.

\cite{BrunelliFedrizzi2024} formulated a property related to the axiomatic systems used for the validation of inconsistency indices based on Pareto dominance.

\cite{Szadoczki2024} examined the efficiency of weight vectors from a geometric viewpoint, and showed (mainly through examples) that the set of efficient weight vectors for a PCM with $3$ alternatives can be represented as a triangle (or a single point in the case of a consistent matrix), while for $4$ alternatives, when there are no consistent triads (triplets) among the alternatives, it can be visualized as a union of three tetrahedra. These tetrahedra are determined by the weight vectors calculated from those incomplete pairwise comparison matrices that can be obtained from the original PCM by deleting some elements, and can be represented with a path spanning tree graph. However, they did not prove the generality of their findings, did not characterize the possible cases, where the geometric properties of the set of efficient weight vectors differ, nor listed the possible cases, only presented different examples. This paper aims to fill in this research gap, provide the first formal proof on the geometry of efficient weight vectors, determine the possible type of geometric characteristic of the set of efficient weight vectors, as well as show the underlying reasons, in the case of $4$ alternatives.

The rest of the paper is organized as follows. Section \ref{sec:2} contains the preliminaries related to pairwise comparison matrices and efficient weight vectors used throughout the study. Section \ref{sec:3} presents the main results, the proof related to the geometric characterization of the set of efficient weight vectors for four alternatives. Finally, Section \ref{sec:4} concludes and raises further research questions. The results of Section \ref{sec:3} are given with more examples and including special cases as well (double perturbed, simple perturbed, consistent matrices) in online Appendices~\ref{append:A} and \ref{append:B}.

\section{Preliminaries}
\label{sec:2}

\begin{definition}[Pairwise comparison matrix (PCM) \cite{Saaty1977}]
    An $\mathbf{A} = \left[ a_{ij} \right]_{i,j=1,\ldots,n}$ square matrix of dimensions $n \times n$ is a pairwise comparison matrix (PCM) if the following properties hold for all $i,j=1,\dots,n$: $a_{ij} > 0$ (positivity) and $a_{ij}=1/a_{ji}$ (reciprocity).
\end{definition}
Note that the second property implies $a_{ii}=1$. One of the most important concept related to pairwise comparison matrices is the consistency.

\begin{definition}[Consistency \cite{Saaty1977}]
    An $\mathbf{A} = \left[ a_{ij} \right]_{i,j=1,\ldots,n}$ PCM is called consistent if $a_{ij}a_{jk}=a_{ik}$ for all $i,j,k=1,\dots,n$.
    A PCM is called inconsistent if it is not consistent.
\end{definition}

There are several methods to obtain a weight vector related to a PCM, among which some of the most popular ones are the right eigenvector method \citep{Saaty1977}, the logarithmic least squares (geometric mean) method \citep{Crawford1985}, the least squares method \citep{Jensen1984}, and the method based on the enumeration of spanning trees \citep{Tsyganok2010,Lundy2017,Mazurek2022} (for further techniques and their analysis, see, for instance \cite{GolanyKress1993,ChooWedley2004,BajwaChooWedley2008,Dijkstra2013}). The ratios of the weights are considered to be the approximations of the elements of the PCM. Different methods can produce different weight vectors. The exact method of obtaining the weight vector will not be important for our purposes.

It may happen that there are missing elements in the PCM. That is called an incomplete pairwise comparison matrix (IPCM).

\begin{definition}[Incomplete pairwise comparison matrix (IPCM)]
An $n\times n$ matrix $\mathbf{A}=[a_{ij}]$ is an incomplete pairwise comparison matrix (IPCM) if:
\begin{itemize}
    \item $a_{ij}\in \mathbb{R}_+ \cup \{\ast\} \ \forall \ 1\leq i,j\leq n$,
    \item $a_{ji}=1/a_{ij}$ if $a_{ij}\in \mathbb{R}_+$,
    \item $a_{ji}=\ast$ if $a_{ij}=\ast$,
\end{itemize}
where $\ast$ denotes the missing elements, and $\mathbb{R}_+$ is the set of positive real numbers.
\end{definition}

Most of the properties of IPCMs can be suitably handled by their representing graphs (graphs of comparisons).

\begin{definition}[Representing graph/Graph of comparisons]
An incomplete pairwise comparison matrix $\mathbf{A}$ can be represented by an undirected graph $G=(V,E)$, where:
\begin{itemize}
    \item the vertices $V=\{1,2,\ldots,n\}$ correspond to the alternatives,
    \item while the edge set $E$ represents the known elements of $\mathbf{A}$ outside the main diagonal:
    $$e_{ij} \in E \iff a_{ij}\neq\ast \ \text{and} \ i\neq j.$$
\end{itemize}
\end{definition}

One can calculate a weight vector from an IPCM as well, where the different methods rely on complementing the matrix as consistently as possible \citep{Tekile2023}. In order to calculate a unique weight vector from an $n\times n$ IPCM, it is necessary that the matrix has a connected representing graph \citep{Bozoki2010}. To reach this, at least $n-1$ comparisons (edges) are needed. Thus, the smallest units from which a unique weight vector can be calculated are the IPCMs that can be represented by a spanning tree (see Example~\ref{representinggraph}), which can always be complemented perfectly consistently. These units also have an important role in determining the set of Pareto efficient weight vectors. 

\begin{example}[Spanning tree representing graph] \label{representinggraph} Let us consider the following incomplete pairwise comparison matrix:

\[
\mathbf{A}=\begin{pmatrix}
1 & a_{12} & \ast  &\ast \\
1/a_{12} & 1 & a_{23} & \ast  \\
\ast  & 1/a_{23} & 1 & a_{34} \\
\ast  & \ast  & 1/a_{34} & 1
\end{pmatrix}
\]

The representing graph $G$ of $\mathbf{A}$ is shown in Figure~\ref{RepresGraph}.

\begin{figure}[H]
    \centering

\begin{tikzpicture}[line join = round, line cap = round]

\tikzstyle{nodeS} = [circle,inner sep=2pt,draw=black,ultra thick];

  \node [nodeS] (B1) at (8,-1) {\footnotesize 1};
\node [nodeS] (B2) at (11,-1) {\footnotesize 2};
\node [nodeS] (B3) at (11,-4) {\footnotesize 3};
\node [nodeS] (B4) at (8,-4) {\footnotesize 4};
  
  \draw [draw=black, very thick](B1) -- (B2);
  \draw [draw=black, very thick](B2) -- (B3);
  \draw [draw=black, very thick](B3) -- (B4);

  \draw (B1) -- (B2) node [midway, above=2pt] {$a_{12}$};
  \draw (B2) -- (B3) node [midway, right=2pt] {$a_{23}$};
  \draw (B3) -- (B4) node [midway, below=2pt] {$a_{34}$};

\end{tikzpicture}

    \caption{The representing graph of $\mathbf{A}$.}
    \label{RepresGraph}
\end{figure}

The exact values of the known comparisons do not matter regarding the representing graph of an IPCM.

\end{example}

Pareto efficiency or Pareto optimality \cite[Chapter 2]{Ehrgott2000}\cite[Chapter 6]{Steuer1986} is an important concept of multi-objective optimization as well as multi-criteria decision making. The Pareto efficiency, also simply called the efficiency of a weight vector corresponding to a PCM can be defined such that by changing the elements of the weight vector, the ratios of the weights cannot approximate any element of the PCM better, without making the approximation of another element worse. The formal definition can be written as follows.

\begin{definition}[Efficient weight vector]
Let $\mathbf{A} =
\left[
a_{ij}
\right]_{i,j=1,\ldots,n}$ be an $n \times n$ PCM and
$\mathbf{w} = (w_1, w_2, \ldots, w_n)^{T}$ be a positive weight vector.
$\mathbf{w}$ is called efficient if no other positive weight vector
$\mathbf{w^{\prime}} = (w^{\prime}_1, w^{\prime}_2, \ldots, w^{\prime}_n)^{T}$
exists such that
\begin{align}
 \left|a_{ij} - \frac{w^{\prime}_i}{w^{\prime}_j} \right| &\leq \left|a_{ij} - \frac{w_i}{w_j} \right| \qquad \text{ for all } 1 \leq i,j \leq n,  \label{eqn:eff1}\\
 \left|a_{k{\ell}} - \frac{w^{\prime}_k}{w^{\prime}_{\ell}} \right| &<  \left|a_{k{\ell}} - \frac{w_k}{w_{\ell}} \right|  \qquad \text{ for some } 1 \leq k,\ell \leq n. \label{eqn:eff2}
\end{align}
A weight vector is called inefficient if it is not efficient.
\end{definition}

From the above definition it follows that renormalization does not affect efficiency or inefficiency, thus we focus on the case, when the weight vector is normalized ($\sum_{i=1}^nw_i=1$).

\begin{remark}
Weight vector $\mathbf{w}$ is efficient if and only if $c\mathbf{w}$ is efficient for any $c > 0$.
\end{remark}

Since for a consistent PCM $\mathbf{A} = \left[ a_{ij} \right]_{i,j=1,\ldots,n}$ there exists (\cite{Saaty1977}) a weight vector $\mathbf{w} = (w_1, w_2, \ldots, w_n)^{T}$ such that $a_{ij}=w_i/w_j$, this $\mathbf{w}$ is the only efficient weight vector (aside from a positive scalar multiplication) in this case.

\cite{Blanquero} devised a necessary and sufficient condition to the efficiency of a weight vector, using the following definition.

\begin{definition}[Blanquero-Carrizosa-Conde (BCC) directed graph]
Let $\mathbf{A}=[a_{ij}]$ be an $n\times n$ PCM, and let $w=\left(w_1,w_2,\ldots,w_n\right)^T$ be a positive weight vector. The Blanquero-Carrizosa-Conde (BCC) directed graph $G=(V,\vv{E})_{\mathbf{A},\mathbf{w}}$ is defined as follows:
\begin{itemize}
    \item the vertices $V=\{1,2,\ldots,n\}$ correspond to the alternatives,
    \item $\vv{E}=\{arc(i\rightarrow j) \ | \ w_i/w_j\geq a_{ij},i\neq j\}$.
\end{itemize}
\end{definition}

The efficiency of a weight vector can then be written by the following necessary and sufficient condition.

\begin{theorem}[{\cite[Corollary 10]{Blanquero}}] \label{thm:TheoremDirectedGraphEfficient} 
Let $\mathbf{A}$ be an $n\times n$ PCM.
A weight vector $\mathbf{w}$ is efficient if and only if the BCC graph
$G=(V,\overrightarrow{E})_{\mathbf{A},\mathbf{w}}$
is a strongly connected digraph, that is,
there exist directed paths from $i$ to $j$ and from $j$ to $i$ for all
pairs of nodes $i , j$.
\end{theorem}

\begin{example}[BCC directed graph]
    \label{BCCexample}
Let us consider the following pairwise comparison matrix $\mathbf{A}$, and weight vector $\mathbf{w}$:

\[
\mathbf{A}=\begin{pmatrix}
1 & a_{12} & a_{13} & a_{14} \\
1/a_{12} & 1 & a_{23} & a_{24} \\
1/a_{13} & 1/a_{23} & 1 & a_{34} \\
1/a_{14} & 1/a_{24} & 1/a_{34} & 1
\end{pmatrix} = 
\begin{pmatrix}
1 & 1 & 5 & 7 \\
1 & 1 & 2 & 8 \\
1/5 & 1/2 & 1 & 1/3 \\
1/7 & 1/8 & 3 & 1
\end{pmatrix},
\]

\[
w=\left(0.25,0.25,0.25,0.25 \right)^T.
\]

The related $G=(V,\overrightarrow{E})_{\mathbf{A},\mathbf{w}}$ BCC directed graph can be seen in Figure~\ref{BCCgraphexample}.

\begin{figure}[H]
    \centering

\begin{tikzpicture}[line join = round, line cap = round]

\tikzstyle{node2} = [circle,inner sep=2pt,draw=black,fill= black,thick];

\node [node2] (A1) at (8.5,-1) {\footnotesize \color{white} 1};
\node [node2] (A2) at (10,-1) {\footnotesize
 \color{white} 2};
 \node [node2] (A3) at (10,-2.5) {\footnotesize \color{white} 3};
\node [node2] (A4) at (8.5,-2.5) {\footnotesize \color{white} 4};

  \draw [latex-latex,draw=black, very thick](A2) -- (A1);
  \draw [-latex,draw=black, very thick](A3) -- (A1);
  \draw [-latex,draw=black, very thick](A3) -- (A2);
  \draw [-latex,draw=black, very thick](A4) -- (A2);
  \draw [-latex,draw=black, very thick](A3) -- (A4);
  \draw [-latex,draw=black, very thick](A4) -- (A1);

\end{tikzpicture}

    \caption{The BCC digraph related to $\mathbf{A}$ and $\mathbf{w}$.}
    \label{BCCgraphexample}
\end{figure}

It can be seen that for the comparison between alternatives $1$ and $2$, the estimation is perfect ($w_1/w_2=a_{12}=1$), thus there are arcs between them with both orientations. However, $G$ is not strongly connected, thus according to Theorem \ref{thm:TheoremDirectedGraphEfficient}, $\mathbf{w}$ is not an efficient weight vector of $\mathbf{A}$.

\end{example}

For any pairs of vertices of the BCC digraph, there is an arc in at least one of the two possible orientations, and there are arcs in both directions only if the comparison of the appropriate alternatives is estimated perfectly by the given weight vector. To determine the set of efficient weight vectors, it is enough to focus on the case when there is exactly one arc for each pair of vertices, i.e., the relevant BCC digraphs are tournaments.

\begin{definition}[Tournament]
    A tournament is a directed graph in which each pair of vertices is connected by an arc with exactly one of the two possible orientations.
\end{definition}

\begin{theorem}[\cite{Camion1959}]
\label{theorem:Hamiltoncycle}
Let $G$ be a tournament. Then, $G$ is strongly connected if and only if $G$ has a Hamiltonian cycle (a cycle that visits each vertex exactly once).
\end{theorem}

In the case when some of the elements are estimated perfectly, there still must be at least one Hamiltonian cycle in the associated BCC digraph to be strongly connected, as every arc has at least one orientation.

\begin{corollary} \label{corolleffhamilton}
    A weight vector $\mathbf{w}$ associated with an $n\times n$ PCM $\mathbf{A}$ is efficient if and only if its BCC digraph $G=(V,\overrightarrow{E})_{\mathbf{A},\mathbf{w}}$ has a Hamiltonian cycle.
\end{corollary}

\begin{remark} \label{remarkHamilton}
For a BCC digraph with $n=3$, there are two possible Hamiltonian cycles (the $3$-cycle with $2$ different orientations), while for $n=4$ there are six possible Hamiltonian cycles (see the cycles in Figure~\ref{Cycles}).
\end{remark}

\section{Results}
\label{sec:3}

The goal of this paper is to determine the possible type of geometric characteristic of the set of efficient weight vectors, and to show the reasons behind that for $4$ alternatives. However, the special cases, when there are consistent cycles in the PCM, are only considered in Appendix~\ref{append:B}, for some geometric interpretation of some of those instances, see also \cite{Szadoczki2024}.

\begin{assumption} \label{assumption:1}
The examined pairwise comparison matrix is a triple perturbed PCM, it has no consistent $k$-cycles ($k\geq3$).
\end{assumption}

A general $4 \times 4$ pairwise comparison matrix has the following form:
\[
\begin{pmatrix}
1 & a_{12} & a_{13} & a_{14} \\
1/a_{12} & 1 & a_{23} & a_{24} \\
1/a_{13} & 1/a_{23} & 1 & a_{34} \\
1/a_{14} & 1/a_{24} & 1/a_{34} & 1
\end{pmatrix},
\]
where $a_{12}, a_{13}, a_{14}, a_{23}, a_{24}, a_{34} >0$. From now on, we use these notations.

As the elements of a positive weight vector $\mathbf{w}$ are usually normalized to one, all of the possible vectors are in the unit $4$-simplex, which can be drawn in $3$ dimensions as a regular tetrahedron (denoted by $S_4$), using the following transformation.

\[ \mathbf{Tw} = \mathbf{\hat{w}}, \] where \[\mathbf{T} =\begin{pmatrix} 1 & 1 & 0 & 0 \\
1 & 0 & 1 & 0 \\
0 & 1 & 1 &  0 \end{pmatrix}.\]

Thus, a given $4$-dimensional weight vector $\mathbf{w}$ can be represented by $\mathbf{\hat{w}}=(x,y,z)^T$ in $S_4$. It is also true that the efficiency of weight vector $\mathbf{w}$ is determined by whether the related BCC graph $G=(V,\overrightarrow{E})_{\mathbf{A},\mathbf{w}}$ is strongly connected. Using the general $4 \times 4$ PCM notations, the arcs of $G$ are determined by the following hyperplanes:
\[h_{12}: \  w_1/w_2 = a_{12}, \]
\[h_{13}: \ w_1/w_3 = a_{13}, \]
\[ \ldots \]
\[h_{34}: \ w_3/w_4 = a_{34}. \]

Figure~\ref{S4} shows $S_4$ with an example plane of $h_{12}:$ $w_1/w_2=a_{12}=2$.

\begin{figure}[H]
    \centering

\begin{tikzpicture}[line join = round, line cap = round, scale=0.82]

\pgfmathsetmacro{\mainfact}{8};
\pgfmathsetmacro{\factor}{1/17};
\pgfmathsetmacro{\factorI}{1/35};
\pgfmathsetmacro{\factorII}{1/13};

\coordinate [label=above left:\text{$V_1=(w_1=1,w_2=0,$}] (A) at (1*\mainfact,1*\mainfact,0);
\coordinate [label=above left:\text{$w_3=0,w_4=0)$}] (AA) at (1.44*\mainfact,1*\mainfact,0);
\coordinate [label=below:\text{$V_2=(w_1=0,w_2=1,w_3=0,w_4=0)$}] (B) at (1*\mainfact,0,1*\mainfact);
\coordinate [label=left:\text{$V_3=(w_1=0,w_2=0,$}] (C) at (0,1*\mainfact,1*\mainfact);
\coordinate [label=left:\text{$w_3=1,w_4=0)$}] (CC) at (0,0.93*\mainfact,1*\mainfact);
\coordinate [label=left:\text{$V_4=(w_1=0,w_2=0,w_3=0,w_4=1)$}] (D) at (0,0,0);

\draw[->] (0,0) -- (10,0,0) node[right] {$x$};
\draw[->] (0,0) -- (0,10,0) node[above] {$y$};
\draw[->] (0,0) -- (0,0,10) node[below left] {$z$};

\foreach \i in {A,B,C,D}
    \draw[dashed,very thick] (0,0)--(\i);
\draw[-,very thick, opacity=.5] (A)--(D)--(B)--cycle;
\draw[-,very thick, opacity=.5] (A) --(D)--(C)--cycle;
\draw[-,very thick, opacity=.5] (B)--(D)--(C)--cycle;
\draw[-,very thick] (A)--(B);
\draw[-,very thick] (A)--(C);
\draw[-,very thick] (A)--(D);

Matrix element 1-2=2
\coordinate [label=right:\text{\color{brown}$h_{12}: \ w_1/w_2=$}] (X1) at (1*\mainfact,2/3*\mainfact,1/3*\mainfact);
\coordinate [label=right:\text{\color{brown}$\hspace{1cm}=a_{12}=2$}] (XX1) at (1*\mainfact,2/3*\mainfact-0.5,1/3*\mainfact);

\draw[-,line width=0.6mm, opacity=.5, fill=brown!30] (C)--(X1)--(D)--cycle;
\draw[-,line width=0.6mm, dashed, draw=brown] (C)--(X1);
\draw[-,line width=0.6mm, draw=brown] (X1)--(D);

\tikzstyle{node1} = [circle,inner sep=2pt,draw=black,fill=black,thick];

\node [node1] (A1) at (2,5) {\footnotesize \color{white} 1};
\node [node1] (A2) at (4,5) {\footnotesize \color{white} 2};

\draw [-latex,draw=brown, very thick](A1) -- (A2);

\node [node1] (A3) at (2.5,1) {\footnotesize \color{white} 1};
\node [node1] (A4) at (4.5,1) {\footnotesize \color{white} 2};

\draw [-latex,draw=brown, very thick](A4) -- (A3);

\end{tikzpicture}

    \caption{$S_4$ with an example plane of $h_{12}:$ $w_1/w_2=a_{12}=2$ and the related arcs of the BCC graph above and below the plane.}
    \label{S4}
\end{figure}

One can see that the possible weight vectors are inside $S_4$ (all the elements of the vectors are positive), while the vertices of $S_4$ ($V_1,V_2,V_3$ and $V_4$) correspond to the cases, where exactly one element is 1, and all the others are 0 (transformation matrix $\mathbf{T}$ is also applied here). The example plane ($h_{12}$) determines the arc of the $G$ BCC graph between alternatives $1$ and $2$, i.e., if $w_1/w_2 > 2$ holds, then there is an arc from $1$ to $2$, and if $w_1/w_2 < 2$, then the arc is from $2$ to $1$. On the plane itself both directions hold. It is also notable that the plane splits the $V_1V_2$ edge of $S_4$ according to a $1:2$ ratio.

For every element of the PCM (aside from the ones coming from the reciprocity property) there is such a plane that determines one of the arcs of the BCC graph. According to Theorem~\ref{thm:TheoremDirectedGraphEfficient}, if the gained BCC graph is strongly connected, then the related weight vector is efficient.

\cite{Szadoczki2024} conjectured via an example that the set of efficient weight vectors for $4$ alternatives geometrically consists of the union of three tetrahedra (denoted by $\mathcal{T}_1,\mathcal{T}_2,$ and $\mathcal{T}_3$) associated with the labeled Hamiltonian cycles of BCC graphs presented in Figure~\ref{Cycles}, where for a given PCM always only one of the two directions is possible for each cycle. This conjecture is proved in the present work.

\begin{figure}[H]
    \centering

\begin{tikzpicture}[line join = round, line cap = round]

\tikzstyle{node1} = [circle,inner sep=2pt,draw=black,fill=black,thick];

\node [node1] (A1) at (4.5,-1) {\footnotesize \color{white} 1};
\node [node1] (A2) at (6,-1) {\footnotesize \color{white}  2};
\node [node1] (A3) at (6,-2.5) {\footnotesize \color{white} 3};
\node [node1] (A4) at (4.5,-2.5) {\footnotesize \color{white} 4};

  \draw [-latex,draw=black, very thick](A1) -- (A2);
  \draw [-latex,draw=black, very thick](A2) -- (A3);
  \draw [-latex,draw=black, very thick](A3) -- (A4);
  \draw [-latex,draw=black, very thick](A4) -- (A1);

\tikzstyle{node2} = [circle,inner sep=2pt,draw=black,fill= black,thick];

\node [node2] (B1) at (8.5,-1) {\footnotesize \color{white} 1};
\node [node2] (B2) at (8.5,-2.5) {\footnotesize \color{white} 3};
\node [node2] (B3) at (10,-2.5) {\footnotesize \color{white} 2};
\node [node2] (B4) at (10,-1) {\footnotesize
 \color{white} 4};

  \draw [-latex,draw=black, very thick](B2) -- (B1);
  \draw [-latex,draw=black, very thick](B3) -- (B2);
  \draw [-latex,draw=black, very thick](B4) -- (B3);
  \draw [-latex,draw=black, very thick](B1) -- (B4);

\tikzstyle{node3} = [circle,inner sep=2pt,draw=black,fill= black,thick];

\node [node3] (C1) at (12.5,-1) {\footnotesize \color{white} 1};
\node [node3] (C2) at (12.5,-2.5) {\footnotesize \color{white} 2};
\node [node3] (C3) at (14,-2.5) {\footnotesize \color{white} 4};
\node [node3] (C4) at (14,-1) {\footnotesize \color{white} 3};

  \draw [-latex,draw=black, very thick](C2) -- (C1);
  \draw [-latex,draw=black, very thick](C3) -- (C2);
  \draw [-latex,draw=black, very thick](C4) -- (C3);
  \draw [-latex,draw=black, very thick](C1) -- (C4);


\node [node1] (AA1) at (4.5,-3.5) {\footnotesize \color{white} 1};
\node [node1] (AA2) at (6,-3.5) {\footnotesize \color{white}  2};
\node [node1] (AA3) at (6,-5) {\footnotesize \color{white} 3};
\node [node1] (AA4) at (4.5,-5) {\footnotesize \color{white} 4};

  \draw [-latex,draw=black, very thick](AA2) -- (AA1);
  \draw [-latex,draw=black, very thick](AA3) -- (AA2);
  \draw [-latex,draw=black, very thick](AA4) -- (AA3);
  \draw [-latex,draw=black, very thick](AA1) -- (AA4);

\tikzstyle{node2} = [circle,inner sep=2pt,draw=black,fill= black,thick];

\node [node2] (BB1) at (8.5,-3.5) {\footnotesize \color{white} 1};
\node [node2] (BB2) at (8.5,-5) {\footnotesize \color{white} 3};
\node [node2] (BB3) at (10,-5) {\footnotesize \color{white} 2};
\node [node2] (BB4) at (10,-3.5) {\footnotesize
 \color{white} 4};

  \draw [-latex,draw=black, very thick](BB1) -- (BB2);
  \draw [-latex,draw=black, very thick](BB2) -- (BB3);
  \draw [-latex,draw=black, very thick](BB3) -- (BB4);
  \draw [-latex,draw=black, very thick](BB4) -- (BB1);

\tikzstyle{node3} = [circle,inner sep=2pt,draw=black,fill= black,thick];

\node [node3] (CC1) at (12.5,-3.5) {\footnotesize \color{white} 1};
\node [node3] (CC2) at (12.5,-5) {\footnotesize \color{white} 2};
\node [node3] (CC3) at (14,-5) {\footnotesize \color{white} 4};
\node [node3] (CC4) at (14,-3.5) {\footnotesize \color{white} 3};

  \draw [-latex,draw=black, very thick](CC1) -- (CC2);
  \draw [-latex,draw=black, very thick](CC2) -- (CC3);
  \draw [-latex,draw=black, very thick](CC3) -- (CC4);
  \draw [-latex,draw=black, very thick](CC4) -- (CC1);


\node (T1) at (5.25,-0) {\footnotesize $\mathcal{T}_1$};
\node (T2) at (9.25,0) {\footnotesize $\mathcal{T}_2$};
\node (T3) at (13.25,0) {\footnotesize $\mathcal{T}_3$};

\end{tikzpicture}

    \caption{The labeled Hamiltonian cycles of BCC graphs associated with tetrahedra $\mathcal{T}_1,\mathcal{T}_2,$ and $\mathcal{T}_3$. Each tetrahedron is related to the cycles below them, from which always only one is possible for a given PCM (see Lemma~\ref{onlyonecycle}).}
    \label{Cycles}
\end{figure}

Each of $\mathcal{T}_1,\mathcal{T}_2,$ and $\mathcal{T}_3$ is determined by exactly four of the six $h_{ij}$ planes ($i=1,\ldots,3, \ j=i+1,\ldots, 4$). 

The question is, which combinations of the Hamiltonian cycles associated with $\mathcal{T}_1,\mathcal{T}_2,$ and $\mathcal{T}_3$ are possible. What are the determining factors behind that, and how can we describe that in a geometric manner? To answer these questions, first the pairwise comparison matrices that should be considered are narrowed down.

It is true that we can get the same geometric object with reversed dimensions (looking at it from different angles), and the properties of pairwise comparisons are often examined through their triads (the different triplets of alternatives) as their building blocks. However, in the examined problem it is not enough to focus on the triads.

\begin{remark} The properties of the set of efficient weight vectors is not solely dependent on the triads of a $4\times4$ pairwise comparison matrix. Furthermore, a PCM can be a double perturbed one (with consistent $4$-cycles violating Assumption~\ref{assumption:1}) even in the case when all of its triads are inconsistent. For further details, see Appendix~\ref{append:A}.

\end{remark}

Instead of the triads, the 4-cycles should be considered, for which the dimension reversion is formalized in Proposition~\ref{negykor}.

\begin{proposition} \label{negykor} 
Any $4\times 4$ pairwise comparison matrix 
\[ \mathbf{A} =
\begin{pmatrix}
1 & a_{12} & a_{13} & a_{14} \\
1/a_{12} & 1 & a_{23} & a_{24} \\
1/a_{13} & 1/a_{23} & 1 & a_{34} \\
1/a_{14} & 1/a_{24} & 1/a_{34} & 1
\end{pmatrix},
\]
for which Assumption~\ref{assumption:1} holds, can be transformed with appropriate reindexing of alternatives to the following form:
\[ \mathbf{A}'' =
\begin{pmatrix}
1 & a_{12}''& a_{13}'' & a_{14}'' \\
1/a_{12}'' & 1 & a_{23}'' & a_{24}'' \\
1/a_{13}'' & 1/a_{23}'' & 1 & a_{34}'' \\
1/a_{14}'' & 1/a_{24}'' & 1/a_{34}'' & 1
\end{pmatrix},
\]
where $\mathbf{A}''$ has the following property:

\bigskip
$a_{12}''a_{23}''a_{34}''a_{41}''<1$, \quad $a_{14}''a_{42}''a_{23}''a_{31}''<1$, \quad  $a_{13}''a_{34}''a_{42}''a_{21}''<1$. \\
where the notation $a_{ij}$ with $i>j$ stands for $1/a_{ji}$.
\end{proposition}
\begin{proof}
\noindent Case 1A: \quad $a_{12}a_{23}a_{34}a_{41}<1$, \quad $a_{14}a_{42}a_{23}a_{31}<1$, \quad  $a_{13}a_{34}a_{42}a_{21}<1$. \\
This case needs no rearrangement, the original $[1,2,3,4]$ satisfies the criteria. Note, that $[1,3,4,2]$ and $[1,4,2,3]$ will also satisfy the criteria in this case, since
\[ \mathbf{A}_{1342} =
\begin{pmatrix}
1 & a_{13} & a_{14} & a_{12} \\
1/a_{13} & 1 & a_{34} & 1/a_{23} \\
1/a_{14} & 1/a_{34} & 1 & 1/a_{24} \\
1/a_{12} & a_{23} & a_{24} & 1
\end{pmatrix}, \
\mathbf{A}_{1423} =
\begin{pmatrix}
1 & a_{14} & a_{12} & a_{13} \\
1/a_{14} & 1 & 1/a_{24} & 1/a_{34} \\
1/a_{12} & a_{24} & 1 & a_{23} \\
1/a_{13} & a_{34} & 1/a_{23} & 1 \\
\end{pmatrix},
\]
and the product of the elements along cycles $(1,2,3,4,1)$, $(1,4,2,3,1)$, and $(1,3,4,2,1)$ are
\[ a_{13} \cdot a_{34} \cdot 1/a_{24} \cdot 1/a_{12} < 1, \quad
a_{12} \cdot a_{23} \cdot a_{34} \cdot 1/a_{14} < 1, \quad
  a_{14} \cdot 1/a_{24} \cdot a_{23} \cdot 1/a_{13} < 1  \]
for matrix $\mathbf{A}_{1342}$ and they are 
\[
a_{14} \cdot 1/a_{24} \cdot a_{23} \cdot 1/a_{13} < 1, \quad
a_{13} \cdot a_{34} \cdot 1/a_{24} \cdot 1/a_{12} < 1 , \quad 
a_{12} \cdot a_{23} \cdot a_{34} \cdot 1/a_{14} < 1
\] for matrix $\mathbf{A}_{1423}$.

The arrangements for this case are when alternative $1$ is put as first, and the directed triad $(2,3,4)$ after that, at any starting point in the latter.

\bigskip
\noindent Case 1B: \quad $a_{12}a_{23}a_{34}a_{41}>1$, \quad $a_{14}a_{42}a_{23}a_{31}>1$, \quad  $a_{13}a_{34}a_{42}a_{21}>1$. \\
Rearrangements $[1,2,4,3]$, $[1,3,2,4]$, and $[1,4,3,2]$ satisfy the requirement.
\[ \mathbf{A}_{1243} =
\begin{pmatrix}
1 & a_{12} & a_{14} & a_{13} \\
1/a_{12} & 1 & a_{24} & a_{23} \\
1/a_{14} & 1/a_{24} & 1 & 1/a_{34}  \\
1/a_{13} & 1/a_{23} & a_{34} & 1
\end{pmatrix}, \
\mathbf{A}_{1324} =
\begin{pmatrix}
1 & a_{13} & a_{12} & a_{14} \\
1/a_{13} & 1 & 1/a_{23} & a_{34} \\
1/a_{12} & a_{23} & 1 & a_{24} \\
1/a_{14} & 1/a_{34} & 1/a_{24} & 1 \\
\end{pmatrix}, \]
\[\mathbf{A}_{1432} =
\begin{pmatrix}
1 & a_{14} & a_{13} & a_{12} \\
1/a_{14} & 1 & 1/a_{34} & 1/a_{24} \\
1/a_{13} & a_{34} & 1 & 1/a_{23} \\
1/a_{12} & a_{24} & a_{23} & 1 \\
\end{pmatrix},
\]
and the product of the elements along cycles $(1,2,3,4,1)$, $(1,4,2,3,1)$ , and $(1,3,4,2,1)$ are
\[ a_{12} \cdot a_{24} \cdot 1/a_{34} \cdot 1/a_{13} < 1, \quad
a_{13} \cdot 1/a_{23} \cdot a_{24} \cdot 1/a_{14} < 1, \quad
a_{14} \cdot 1/a_{34} \cdot 1/a_{23} \cdot 1/a_{12} < 1 \]
for matrix $\mathbf{A}_{1243}$, while they are
\[ a_{13} \cdot 1/a_{23} \cdot a_{24} \cdot 1/a_{14} < 1, \quad
a_{14} \cdot 1/a_{34} \cdot 1/a_{23} \cdot 1/a_{12} < 1, \quad
a_{12} \cdot a_{24} \cdot 1/a_{34} \cdot 1/a_{13} < 1 \]
for matrix $\mathbf{A}_{1324}$ and they are
\[ a_{14} \cdot 1/a_{34} \cdot 1/a_{23} \cdot 1/a_{12} < 1, \quad
a_{12} \cdot a_{24} \cdot 1/a_{34} \cdot 1/a_{13} < 1, \quad
a_{13} \cdot 1/a_{23} \cdot a_{24} \cdot 1/a_{14} < 1
\]
for matrix $\mathbf{A}_{1432}$.

The arrangements for this case are when alternative $1$ is put as first and the directed triad $(4,3,2)$ after that, at any starting point in the latter.

\bigskip
\noindent Case 2A: \quad $a_{12}a_{23}a_{34}a_{41}<1$, \quad $a_{14}a_{42}a_{23}a_{31}>1$, \quad  $a_{13}a_{34}a_{42}a_{21}<1$. \\
Rearrangements $[2,1,3,4]$, $[2,3,4,1]$, and $[2,4,1,3]$ satisfy the requirement.
\[ \mathbf{A}_{2134} =
\begin{pmatrix}
1 & 1/a_{12} & a_{23} & a_{24} \\
a_{12} & 1 & a_{13} & a_{14} \\
1/a_{23} & 1/a_{13} & 1 & a_{34} \\
1/a_{24} & 1/a_{14} & 1/a_{34} & 1 \\
\end{pmatrix}, \
\mathbf{A}_{2341} =
\begin{pmatrix}
1 & a_{23} & a_{24} & 1/a_{12} \\
1/a_{23} & 1 & a_{34} & 1/a_{13} \\
1/a_{24} & 1/a_{34} & 1 & 1/a_{14} \\
a_{12} & a_{13} & a_{14} & 1 \\
\end{pmatrix}, \]
\[\mathbf{A}_{2413} =
\begin{pmatrix}
1 & a_{24} & 1/a_{12} & a_{23} \\
1/a_{24} & 1 & 1/a_{14} & 1/a_{34} \\
a_{12} & a_{14} & 1 & a_{13} \\
1/a_{23} & a_{34} & 1/a_{13} & 1 \\
\end{pmatrix},
\]
and the product of the elements along cycles $(1,2,3,4,1)$, $(1,4,2,3,1)$, and $(1,3,4,2,1)$ are
\[  1/a_{12} \cdot a_{13} \cdot a_{34} \cdot 1/a_{24} < 1, \quad 
 a_{24} \cdot 1/a_{14} \cdot a_{13} \cdot 1/a_{23} < 1,  \quad
a_{23} \cdot a_{34} \cdot 1/a_{14} \cdot a_{12} < 1
\]
for matrix $\mathbf{A}_{2134}$, while they are
\[ a_{23} \cdot a_{34} \cdot 1/a_{14} \cdot a_{12} < 1, \quad
 1/a_{12} \cdot a_{13} \cdot a_{34} \cdot 1/a_{24} < 1, \quad
 a_{24} \cdot 1/a_{14} \cdot a_{13} \cdot 1/a_{23} < 1
 \]
for matrix $\mathbf{A}_{2341}$ and they are
\[ a_{24} \cdot 1/a_{14} \cdot a_{13} \cdot 1/a_{23} < 1, \quad
a_{23} \cdot a_{34} \cdot 1/a_{14} \cdot a_{12} < 1, \quad
1/a_{12} \cdot a_{13} \cdot a_{34} \cdot 1/a_{24} < 1
\]
for matrix $\mathbf{A}_{2413}$.

The arrangements for this case are when alternative $2$ is put as first and the directed triad $(1,3,4)$ after that, at any starting point in the latter.

\bigskip
\noindent Case 2B: \quad $a_{12}a_{23}a_{34}a_{41}>1$, \quad $a_{14}a_{42}a_{23}a_{31}<1$, \quad  $a_{13}a_{34}a_{42}a_{21}>1$. \\
Rearrangements $[2,1,4,3]$, $[2,4,3,1]$, and $[2,3,1,4]$ satisfy the requirement.
\[ \mathbf{A}_{2143} =
\begin{pmatrix}
1 & 1/a_{12} & a_{24} & a_{23} \\
a_{12} & 1 & a_{14} & a_{13} \\
1/a_{24} & 1/a_{14} & 1 & 1/a_{34} \\
1/a_{23} & 1/a_{13} & a_{34} & 1 \\
\end{pmatrix}, \
\mathbf{A}_{2431} =
\begin{pmatrix}
1 & a_{24} & a_{23} & 1/a_{12} \\
1/a_{24} & 1 & 1/a_{34} & 1/a_{14} \\
1/a_{23} & a_{34} & 1 & 1/a_{13} \\
a_{12} & a_{14} & a_{13} & 1 \\
\end{pmatrix}, \]
\[\mathbf{A}_{2314} =
\begin{pmatrix}
1 & a_{23} & 1/a_{12} & a_{24} \\
1/a_{23} & 1 & 1/a_{13} & a_{34} \\
a_{12} & a_{13} & 1 & a_{14} \\
1/a_{24} & 1/a_{34} & 1/a_{14} & 1 \\
\end{pmatrix},
\]
and the product of the elements along cycles $(1,2,3,4,1)$, $(1,4,2,3,1)$, and $(1,3,4,2,1)$ are
\[ 1/a_{12} \cdot a_{14} \cdot 1/a_{34} \cdot 1/a_{23} < 1, \quad
a_{23} \cdot 1/a_{13} \cdot a_{14} \cdot 1/a_{24} < 1, \quad
a_{24} \cdot 1/a_{34} \cdot 1/a_{13} \cdot a_{12} < 1 \]
for matrix $\mathbf{A}_{2143}$, while they are
\[ a_{24} \cdot 1/a_{34} \cdot 1/a_{13} \cdot a_{12} < 1, \quad
1/a_{12} \cdot a_{14} \cdot 1/a_{34} \cdot 1/a_{23} < 1, \quad
a_{23} \cdot 1/a_{13} \cdot a_{14} \cdot 1/a_{24} < 1 \]
for matrix $\mathbf{A}_{2431}$ and they are
\[ a_{23} \cdot 1/a_{13} \cdot a_{14} \cdot 1/a_{24} < 1, \quad
a_{24} \cdot 1/a_{34} \cdot 1/a_{13} \cdot a_{12} < 1, \quad
1/a_{12} \cdot a_{14} \cdot 1/a_{34} \cdot 1/a_{23} < 1
\]
for matrix $\mathbf{A}_{2314}$.

The arrangements for this case are when alternative $2$ is put as first and the directed triad $(4,3,1)$ after that, at any starting point in the latter.

\bigskip
\noindent Case 3A: \quad $a_{12}a_{23}a_{34}a_{41}>1$, \quad $a_{14}a_{42}a_{23}a_{31}<1$, \quad  $a_{13}a_{34}a_{42}a_{21}<1$. \\
Rearrangements $[3,1,4,2]$, $[3,4,2,1]$, and $[3,2,1,4]$ satisfy the requirement.
\[ \mathbf{A}_{3142} =
\begin{pmatrix}
1 & 1/a_{13} & a_{34} & 1/a_{23} \\
a_{13} & 1 & a_{14} & a_{12} \\
1/a_{34} & 1/a_{14} & 1 & 1/a_{24} \\
a_{23} & 1/a_{12} & a_{24} & 1 \\
\end{pmatrix}, \
\mathbf{A}_{3421} =
\begin{pmatrix}
1 & a_{34} & 1/a_{23} & 1/a_{13} \\
1/a_{34} & 1 & 1/a_{24} & 1/a_{14} \\
a_{23} & a_{24} & 1 & 1/a_{12} \\
a_{13} & a_{14} & a_{12} & 1 \\
\end{pmatrix}, \]
\[\mathbf{A}_{3214} =
\begin{pmatrix}
1 & 1/a_{23} & 1/a_{13} & a_{34} \\
a_{23} & 1 & 1/a_{12} & a_{24} \\
a_{13} & a_{12} & 1 & a_{14} \\
1/a_{34} & 1/a_{24} & 1/a_{14} & 1 \\
\end{pmatrix},
\]
and the product of the elements along cycles $(1,2,3,4,1)$, $(1,4,2,3,1)$, and $(1,3,4,2,1)$ are
\[ 1/a_{13} \cdot a_{14} \cdot 1/a_{24} \cdot a_{23} < 1, \quad
1/a_{23} \cdot 1/a_{12} \cdot a_{14} \cdot 1/a_{34} < 1, \quad
a_{34} \cdot 1/a_{24} \cdot 1/a_{12} \cdot a_{13} < 1 \]
for matrix $\mathbf{A}_{3142}$, while they are
\[ a_{34} \cdot 1/a_{24} \cdot 1/a_{12} \cdot a_{13} < 1, \quad
1/a_{13} \cdot a_{14} \cdot 1/a_{24} \cdot a_{23} < 1, \quad
1/a_{23} \cdot 1/a_{12} \cdot a_{14} \cdot 1/a_{34} < 1 \]
for matrix $\mathbf{A}_{3421}$ and they are
\[ 1/a_{23} \cdot 1/a_{12} \cdot a_{14} \cdot 1/a_{34} < 1, \quad
a_{34} \cdot 1/a_{24} \cdot 1/a_{12} \cdot a_{13} < 1, \quad
1/a_{13} \cdot a_{14} \cdot 1/a_{24} \cdot a_{23} < 1 \]
for matrix $\mathbf{A}_{3214}$.

The arrangements for this case are when alternative $3$ is put as first and the directed triad $(4,2,1)$ after that, at any starting point in the latter.

\bigskip
\noindent Case 3B: \quad $a_{12}a_{23}a_{34}a_{41}<1$, \quad $a_{14}a_{42}a_{23}a_{31}>1$, \quad  $a_{13}a_{34}a_{42}a_{21}>1$. \\
Rearrangements $[3,1,2,4]$, $[3,2,4,1]$, and $[3,4,1,2]$ satisfy the requirement.
\[ \mathbf{A}_{3124} =
\begin{pmatrix}
1 & 1/a_{13} & 1/a_{23} & a_{34} \\
a_{13} & 1 & a_{12} & a_{14} \\
a_{23} & 1/a_{12} & 1 & a_{24} \\
1/a_{34} & 1/a_{14} & 1/a_{24} & 1 \\
\end{pmatrix}, \
\mathbf{A}_{3241} =
\begin{pmatrix}
1 & 1/a_{23} & a_{34} & 1/a_{13} \\
a_{23} & 1 & a_{24} & 1/a_{12} \\
1/a_{34} & 1/a_{24} & 1 & 1/a_{14} \\
a_{13} & a_{12} & a_{14} & 1 \\
\end{pmatrix}, \]
\[\mathbf{A}_{3412} =
\begin{pmatrix}
1 & a_{34} & 1/a_{13} & 1/a_{23} \\
1/a_{34} & 1 & 1/a_{14} & 1/a_{24} \\
a_{13} & a_{14} & 1 & a_{12} \\
a_{23} & a_{24} & 1/a_{12} & 1 \\
\end{pmatrix},
\]
and the product of the elements along cycles $(1,2,3,4,1)$, $(1,4,2,3,1)$, and $(1,3,4,2,1)$ are
\[ 1/a_{13} \cdot a_{12} \cdot a_{24} \cdot 1/a_{34} < 1, \quad
a_{34} \cdot 1/a_{14} \cdot a_{12} \cdot a_{23} < 1, \quad
1/a_{23} \cdot a_{24} \cdot 1/a_{14} \cdot a_{13} < 1 \]
for matrix $\mathbf{A}_{3124}$, while they are
\[ 1/a_{23} \cdot a_{24} \cdot 1/a_{14} \cdot a_{13} < 1, \quad
1/a_{13} \cdot a_{12} \cdot a_{24} \cdot 1/a_{34} < 1, \quad
a_{34} \cdot 1/a_{14} \cdot a_{12} \cdot a_{23} < 1 \]
for matrix $\mathbf{A}_{3241}$ and they are
\[ a_{34} \cdot 1/a_{14} \cdot a_{12} \cdot a_{23} < 1, \quad
1/a_{23} \cdot a_{24} \cdot 1/a_{14} \cdot a_{13} < 1, \quad
1/a_{13} \cdot a_{12} \cdot a_{24} \cdot 1/a_{34} < 1 \]
for matrix $\mathbf{A}_{3412}$.

The arrangements for this case are when alternative $3$ is put as first and the directed triad $(1,2,4)$ after that, at any starting point in the latter.

\bigskip
\noindent Case 4A: \quad $a_{12}a_{23}a_{34}a_{41}<1$, \quad $a_{14}a_{42}a_{23}a_{31}<1$, \quad  $a_{13}a_{34}a_{42}a_{21}>1$. \\
Rearrangements $[4,1,2,3]$, $[4,2,3,1]$, and $[4,3,1,2]$ satisfy the requirement.
\[ \mathbf{A}_{4123} =
\begin{pmatrix}
1 & 1/a_{14} & 1/a_{24} & 1/a_{34} \\
a_{14} & 1 & a_{12} & a_{13} \\
a_{24} & 1/a_{12} & 1 & a_{23} \\
a_{34} & 1/a_{13} & 1/a_{23} & 1 \\
\end{pmatrix}, \
\mathbf{A}_{4231} =
\begin{pmatrix}
1 & 1/a_{24} & 1/a_{34} & 1/a_{14} \\
a_{24} & 1 & a_{23} & 1/a_{12} \\
a_{34} & 1/a_{23} & 1 & 1/a_{13} \\
a_{14} & a_{12} & a_{13} & 1 \\
\end{pmatrix}, \]
\[\mathbf{A}_{4312} =
\begin{pmatrix}
1 & 1/a_{34} & 1/a_{14} & 1/a_{24} \\
a_{34} & 1 & 1/a_{13} & 1/a_{23} \\
a_{14} & a_{13} & 1 & a_{12} \\
a_{24} & a_{23} & 1/a_{12} & 1 \\
\end{pmatrix},
\]
and the product of the elements along cycles $(1,2,3,4,1)$, $(1,4,2,3,1)$, and $(1,3,4,2,1)$ are
\[ 1/a_{14} \cdot a_{12} \cdot a_{23} \cdot a_{34} < 1, \quad
1/a_{34} \cdot 1/a_{13} \cdot a_{12} \cdot a_{24} < 1, \quad
1/a_{24} \cdot a_{23} \cdot 1/a_{13} \cdot a_{14} < 1 \]
for matrix $\mathbf{A}_{4123}$, while they are
\[1/a_{24} \cdot a_{23} \cdot 1/a_{13} \cdot a_{14} < 1, \quad
1/a_{14} \cdot a_{12} \cdot a_{23} \cdot a_{34} < 1,\quad
1/a_{34} \cdot 1/a_{13} \cdot a_{12} \cdot a_{24} < 1 \]
for matrix $\mathbf{A}_{4231}$ and they are
\[ 1/a_{34} \cdot 1/a_{13} \cdot a_{12} \cdot a_{24} < 1, \quad
1/a_{24} \cdot a_{23} \cdot 1/a_{13} \cdot a_{14} < 1, \quad
1/a_{14} \cdot a_{12} \cdot a_{23} \cdot a_{34} < 1 \]
for matrix $\mathbf{A}_{4312}$.

The arrangements for this case are when alternative $4$ is put as first and the directed triad $(1,2,3)$ after that, at any starting point in the latter.

\bigskip
\noindent Case 4B: \quad $a_{12}a_{23}a_{34}a_{41}>1$, \quad $a_{14}a_{42}a_{23}a_{31}>1$, \quad  $a_{13}a_{34}a_{42}a_{21}<1$. \\
Rearrangements $[4,1,3,2]$, $[4,2,1,3]$, and $[4,3,2,1]$ satisfy the requirement.
\[ \mathbf{A}_{4132} =
\begin{pmatrix}
1 & 1/a_{14} & 1/a_{34} & 1/a_{24} \\
a_{14} & 1 & a_{13} & a_{12} \\
a_{34} & 1/a_{13} & 1 & 1/a_{23} \\
a_{24} & 1/a_{12} & a_{23} & 1 \\
\end{pmatrix}, \
\mathbf{A}_{4213} =
\begin{pmatrix}
1 & 1/a_{24} & 1/a_{14} & 1/a_{34} \\
a_{24} & 1 & 1/a_{12} & a_{23} \\
a_{14} & a_{12} & 1 & a_{13} \\
a_{34} & 1/a_{23} & 1/a_{13} & 1 \\
\end{pmatrix}, \]
\[\mathbf{A}_{4321} =
\begin{pmatrix}
1 & 1/a_{34} & 1/a_{24} & 1/a_{14} \\
a_{34} & 1 & 1/a_{23} & 1/a_{13} \\
a_{24} & a_{23} & 1 & 1/a_{12} \\
a_{14} & a_{13} & a_{12} & 1 \\
\end{pmatrix},
\]
and the product of the elements along cycles $(1,2,3,4,1)$, $(1,4,2,3,1)$, and $(1,3,4,2,1)$ are
\[1/a_{14} \cdot a_{13} \cdot 1/a_{23} \cdot a_{24} < 1, \quad
1/a_{24} \cdot 1/a_{12} \cdot a_{13} \cdot a_{34} < 1 , \quad
1/a_{34} \cdot 1/a_{23} \cdot 1/a_{12} \cdot a_{14} < 1 \]
for matrix $\mathbf{A}_{4132}$, while they are
\[1/a_{24} \cdot 1/a_{12} \cdot a_{13} \cdot a_{34} < 1, \quad
1/a_{34} \cdot 1/a_{23} \cdot 1/a_{12} \cdot a_{14} < 1,\quad
1/a_{14} \cdot a_{13} \cdot 1/a_{23} \cdot a_{24} < 1 \]
for matrix $\mathbf{A}_{4213}$ and they are
\[ 1/a_{34} \cdot 1/a_{23} \cdot 1/a_{12} \cdot a_{14} < 1, \quad
1/a_{14} \cdot a_{13} \cdot 1/a_{23} \cdot a_{24} < 1, \quad
1/a_{24} \cdot 1/a_{12} \cdot a_{13} \cdot a_{34} < 1 \]
for matrix $\mathbf{A}_{4321}$.

The arrangements for this case are when alternative $4$ is put as first and the directed triad $(3,2,1)$ after that, at any starting point in the latter.
\end{proof}

That means any PCM can be rearranged in a way that the directed $4$-cycles have a given orientation (see Corollary~\ref{A''orientations}). In order to prove that the set of efficient weight vectors is determined by three tetrahedra ($\mathcal{T}_1$, $ \mathcal{T}_2$, and $\mathcal{T}_3$) a number of lemmas will be necessary. 

\begin{lemma}  \label{spanningtrees}
For a $4\times4$ pairwise comparison matrix $\mathbf{A}$ the weight vectors calculated from any of its incomplete submatrices which can be represented by a spanning tree will be efficient, and each will estimate three elements of $\mathbf{A}$ perfectly.
\end{lemma}

\begin{proof}

An incomplete pairwise comparison submatrix of $\mathbf{A}$ that can be represented by a spanning tree contains three comparisons and there are no $k$-cycles ($k\geq3$) in it. That means there cannot be inconsistency between them, the given submatrix can be complemented in a consistent way. Hence the $\mathbf{w}$ weight vector calculated from it will estimate the three given elements perfectly.

There cannot be a $\mathbf{w}'$ that estimates these three elements perfectly and at least one other element better than $\mathbf{w}$, as the single normalized weight vector that can estimate three elements of the matrix (and thus the consistent PCM calculated from them) perfectly is $\mathbf{w}$. Thus, $\mathbf{w}$ is efficient for $\mathbf{A}$.
\end{proof}

\begin{lemma} \label{onlyonecycle}
    For the BCC graphs of a $4\times4$ pairwise comparison matrix $\mathbf{A}$ (that satisfies Assumption~\ref{assumption:1}), from the six labeled Hamiltonian directed cycles exactly three are possible, which are the three labeled Hamiltonian undirected cycles with exactly one orientation each.
\end{lemma}

\begin{proof}
    There are three labeled Hamiltonian undirected cycles for four alternatives: $(1,2,3,4,1)$, $(1,4,2,3,1)$, and $(1,3,4,2,1)$. Let us consider the $(1,2,3,4,1)$ cycle, which has two possible orientations.

    Without the loss of generality, let us assume that there exists a (positive) weight vector $\mathbf{w}$ for which the orientation of the $G=(V,\overrightarrow{E})_{\mathbf{A},\mathbf{w}}$ BCC graph is $(1,2,3,4,1)$. That means the following:
    \[w_1/w_2\geq a_{12} \iff w_1\geq a_{12}w_2 \tag{I} \label{I}\]
    \[w_2/w_3\geq a_{23} \iff w_2\geq a_{23}w_3 \tag{II} \label{II}\]
    \[w_3/w_4\geq a_{34} \iff w_3 \geq a_{34}w_4 \tag{III} \label{III}\]
    \[w_4/w_1\geq a_{41} \iff w_4\geq a_{41}w_1 \tag{IV} \label{IV}\]

    Writing Inequality~\ref{II} into \ref{I}, and then \ref{III} into the result, and finally \ref{IV} into the result, we get:
    \[w_1\geq a_{12}a_{23}a_{34}a_{41}w_1 \iff 1\geq a_{12}a_{23}a_{34}a_{41} \tag{V} \label{V}\]

    If we assume a weight vector with the other orientation of the BCC graph, then Inequality~\ref{V} will be reversed. The only possible way to get BCC graphs with both orientations is to get an equality in \ref{V}. However, that means the $(1,2,3,4,1)$ $4$-cycle is consistent, and that violates Assumption~\ref{assumption:1}.

    The same logic can be used for the other two $4$-cycles.
\end{proof}

\begin{corollary} \label{A''orientations} For any $4\times 4$ pairwise comparison matrix $\mathbf{A}$ (that satisfies Assumption~\ref{assumption:1}), for the $\mathbf{A}''$ matrix obtained from $\mathbf{A}$ by appropriate rearrangements shown in Proposition~\ref{negykor}, the three possible directed Hamiltonian cycles are $(1,2,3,4,1)$, $(1,4,2,3,1)$, and $(1,3,4,2,1)$.
    
\end{corollary}

\begin{proposition} \label{final} For any $4\times4$ pairwise comparison matrix $\mathbf{A}$ (that satisfies Assumption~\ref{assumption:1}), there exists $\mathbf{A''}$ obtained from $\mathbf{A}$ by appropriate rearrangements, for which the set of efficient weight vectors can be constructed as the union of three tetrahedra ($\mathcal{T}_1,\mathcal{T}_2,$ and $\mathcal{T}_3$), where each tetrahedron is determined by four weight vectors, each calculated from an incomplete pairwise comparison submatrix of  $\mathbf{A''}$ with a path spanning tree representing graph. Additionally, for each efficient weight vector at least one of the following directed cycles is present in its BCC graph: (1,2,3,4,1), (1,4,2,3,1), or (1,3,4,2,1).
\end{proposition}

\begin{proof}
Based on Corollary~\ref{corolleffhamilton} the BCC graph of any efficient weight vector must contain a directed Hamiltonian cycle. According to Lemma~\ref{onlyonecycle} for a given $4\times4$ PCM $\mathbf{A}$, exactly three directed Hamiltonian cycles are possible. Based on Lemma~\ref{spanningtrees} the weight vectors calculated from an incomplete submatrix of $\mathbf{A}$ with a spanning tree representing graph must be efficient, and

\begin{itemize}
    \item each estimates three elements of $\mathbf{A}$ perfectly
    
    $\Longupdownarrow$
    
    \item they are at the intersection of three of the possible six cutting planes $h_{12},\ldots,h_{34}$ 

$\Longupdownarrow$
    
    \item in their BCC graphs there are arcs in both directions between three pairs of vertices.  
\end{itemize}
Apply the reindexing of Proposition~\ref{negykor} to $\mathbf{A}$ to gain $\mathbf{A}''$, where the three possible directed Hamiltonian cycles' orientations are: $(1,2,3,4,1)$, $(1,4,2,3,1)$ and $(1,3,4,2,1)$ as stated in Corollary~\ref{A''orientations}. 

Without the loss of generality, let us consider the $(1,2,3,4,1)$ cycle (that is associated with $\mathcal{T}_1$). There are four incomplete pairwise comparison submatrices of $\mathbf{A}''$ that have path spanning tree representing graphs (denoted by $T_1$, $T_2$, $T_3$, and $T_4$ respectively) that include only comparisons related to this cycle:
    \begin{itemize}
        \item $T_1: (1,2,3,4)$, excluding the comparison of $(1,4)$,
        \item $T_2: (2,3,4,1)$, excluding the comparison of $(1,2)$,
        \item $T_3: (3,4,1,2)$, excluding the comparison of $(2,3)$,
        \item $T_4: (4,1,2,3)$, excluding the comparison of $(3,4)$.
    \end{itemize}

Based on Assumption~\ref{assumption:1} there are no consistent $k$-cycles ($k\geq3$) in $\mathbf{A}''$, thus the weight vectors calculated from $T_1$, $T_2$, $T_3$, and $T_4$ (denoted by $\mathbf{w}^{T_1}$, $\mathbf{w}^{T_2}$, $\mathbf{w}^{T_3}$, and $\mathbf{w}^{T_4}$ respectively) do not correspond to each other. They are efficient, and as four different points, they determine a tetrahedron ($\mathcal{T}_1$) in the space of normalized ($\sum_{i=1}^4w_i=1$) $4$-dimensional weight vectors.

 We would like to show that a given $\mathbf{w}$ weight vector generates a BCC graph that contains the $(1,2,3,4,1)$ directed Hamiltonian cycle if and only if it is in tetrahedron $\mathcal{T}_1$ determined by $\mathbf{w}^{T_1}$, $\mathbf{w}^{T_2}$, $\mathbf{w}^{T_3}$, and $\mathbf{w}^{T_4}$, that is, $\mathbf{w}$ is in their convex hull.

First, let us consider the nodes of the tetrahedron ($\mathbf{w}^{T_1}$, $\mathbf{w}^{T_2}$, $\mathbf{w}^{T_3}$, and $\mathbf{w}^{T_4}$). Each of them estimates three elements of $\mathbf{A}''$ perfectly, which correspond to three edges along the $(1,2,3,4,1)$ cycle, which means that the corresponding edges of their BCC graphs have both possible orientations. However, whichever orientation does the last edge has in this cycle it will result in either the $(1,2,3,4,1)$ oriented cycle, or the opposite $(1,4,3,2,1)$ cycle. However, based on Lemma~\ref{onlyonecycle} only one of those is possible, which is $(1,2,3,4,1)$ for $\mathbf{A}''$. Thus, all of $\mathbf{w}^{T_1}$, $\mathbf{w}^{T_2}$, $\mathbf{w}^{T_3}$, and $\mathbf{w}^{T_4}$ generate BCC graphs that contain the $(1,2,3,4,1)$ cycle, they also satisfy the relations described by (\ref{i}).
\[\mathbf{w}^{T_1}: \hspace{2.8cm} \mathbf{w}^{T_2}: \hspace{2.8cm} \mathbf{w}^{T_3}: \hspace{2.8cm} \mathbf{w}^{T_4}:\]
\[w^{T_1}_1/w^{T_1}_2=a_{12}'' \quad \quad \quad w^{T_2}_1/w^{T_2}_2>a_{12}'' \quad \quad \quad w^{T_3}_1/w^{T_3}_2=a_{12}'' \quad \quad \quad w^{T_4}_1/w^{T_4}_2=a_{12}''\]
\[w^{T_1}_2/w^{T_1}_3=a_{23}'' \quad \quad \quad w^{T_2}_2/w^{T_2}_3=a_{23}'' \quad \quad \quad w^{T_3}_2/w^{T_3}_3>a_{23}'' \quad \quad \quad w^{T_4}_2/w^{T_4}_3=a_{23}'' \tag{i} \label{i}\]
\[w^{T_1}_3/w^{T_1}_4=a_{34}'' \quad \quad \quad w^{T_2}_3/w^{T_2}_4=a_{34}'' \quad \quad \quad w^{T_3}_3/w^{T_3}_4=a_{34}'' \quad \quad \quad w^{T_4}_3/w^{T_4}_4>a_{34}''\]
\[w^{T_1}_4/w^{T_1}_1>a_{41}'' \quad \quad \quad w^{T_2}_4/w^{T_2}_1=a_{41}'' \quad \quad \quad w^{T_3}_4/w^{T_3}_1=a_{41}'' \quad \quad \quad w^{T_4}_4/w^{T_4}_1=a_{41}''\]

The BCC graph generated by the $\mathbf{w}$ weight vector contains the $(1,2,3,4,1)$ directed cycle if the following holds:
\[w_1/w_2\geq a_{12}'',\]
\[w_2/w_3\geq a_{23}'', \tag{ii} \label{ii}\]
\[w_3/w_4\geq a_{34}'',\]
\[w_4/w_1\geq a_{41}''.\]

Based on (\ref{i}) any $\mathbf{w}$ weight vector that is in the tetrahedron determined by $\mathbf{w}^{T_1}$, $\mathbf{w}^{T_2}$, $\mathbf{w}^{T_3}$, and $\mathbf{w}^{T_4}$ generates a BCC graph containing the $(1,2,3,4,1)$ directed cycle (it satisfies (\ref{ii})), as there exist nonnegative $\lambda_1$, $\lambda_2$, $\lambda_3$, and $\lambda_4$ ($\sum_{i=1}^4 \lambda_i=1$) for which:
\[\mathbf{w}=\lambda_1\mathbf{w}^{T_1}+\lambda_2\mathbf{w}^{T_2}+\lambda_3\mathbf{w}^{T_3}+\lambda_4\mathbf{w}^{T_4}.\]

Each of the faces of the tetrahedron are determined by three of $\mathbf{w}^{T_1}$, $\mathbf{w}^{T_2}$, $\mathbf{w}^{T_3}$, and $\mathbf{w}^{T_4}$. There is exactly one of the four types of inequalities in (\ref{i}) that will be an equation for all three that determines a given face of the tetrahedron. However, that means the given face of the tetrahedron is on a cutting plane that determines the orientation of a given edge related to the $(1,2,3,4,1)$ cycle. We just showed that inside the tetrahedron this cycle holds, thus the given edge must be reversed outside of the tetrahedron. Thus, it is also true that any $\mathbf{w}$ that generates a BCC graph containing the $(1,2,3,4,1)$ directed cycle must be in this tetrahedron.

The same logic can be applied to the directed cycles $(1,4,2,3,1)$ (associated with $\mathcal{T}_2$) and $(1,3,4,2,1)$ (associated with $\mathcal{T}_3$).

\end{proof}

\begin{corollary} \label{convexcomb} For any $4\times4$ pairwise comparison matrix $\mathbf{A}$ (that satisfies Assumption~\ref{assumption:1}), every efficient weight vector can be written as the convex combination of the $4$ weight vectors determined by the path spanning tree representing graphs for one of the three cycles of  (1,2,3,4,1), (1,4,2,3,1), or (1,3,4,2,1).
    
\end{corollary}

\begin{remark} \label{connection} For any $4\times4$ pairwise comparison matrix $\mathbf{A}$ (that satisfies Assumption~\ref{assumption:1}) for each pair of $\mathcal{T}_1$, $\mathcal{T}_2$ and $\mathcal{T}_3$, for the two tetrahedra in the given pair, one of their edges is on the same line, and two of their faces are pairwise on the same planes. For further details, see Proposition~\ref{setofefficientweightvectorsspeccases} in Appendix~\ref{append:B}.
    
\end{remark}

\begin{remark}\label{speccases} The results of Proposition~\ref{negykor}, Proposition~\ref{final}, and Corollary~\ref{convexcomb} hold without Assumption~\ref{assumption:1} with appropriate modifications (allowing equaility as well). Thus, the main results are the same even if there are consistent cycles in the PCM, however, each consistent subset will  only result in one possible efficient weight vector. For all the special cases with examples, see Appendix~\ref{append:B}.

\end{remark}

Let us consider the set of efficient weight vectors geometrically for the pairwise comparison matrix used in Example~\ref{BCCexample}.

\begin{example} \label{mainexample}
    
As mentioned before, each of $\mathcal{T}_1,\mathcal{T}_2,$ and $\mathcal{T}_3$ is determined by exactly four cutting planes. Figures~\ref{Cycle1Final}, \ref{Cycle2Final}, and \ref{Cycle3Final} show $\mathcal{T}_1,\mathcal{T}_2,$ and $\mathcal{T}_3$, respectively, with the planes determining them, while Figure~\ref{OnlyOptFinal} presents the set of efficient weight vectors (the union of $\mathcal{T}_1,\mathcal{T}_2,$ and $\mathcal{T}_3$) for the following illustrative example matrix ($a_{12}'' = 1$, $a_{13}''=5$, $a_{14}''=7$, $a_{23}''=2$, $a_{24}'' = 8$, and $a_{34}''=1/3$):

\[ \mathbf{A''} =
\begin{pmatrix}
1 & 1 & 5 & 7 \\
1 & 1 & 2 & 8 \\
1/5 & 1/2 & 1 & 1/3 \\
1/7 & 1/8 & 3 & 1
\end{pmatrix}.
\]

\begin{figure}[H]
    \centering

\begin{tikzpicture}[line join = round, line cap = round, scale=0.87]

\pgfmathsetmacro{\mainfact}{8};

\coordinate [label=above left:\text{$V_1=(w_1=1,w_2=0,w_3=0,w_4=0)$}] (A) at (1*\mainfact,1*\mainfact,0);
\coordinate [label=below:\text{$V_2=(w_1=0,w_2=1,w_3=0,w_4=0)$}] (B) at (1*\mainfact,0,1*\mainfact);
\coordinate [label=left:\text{$V_3=(w_1=0,w_2=0,$}] (C) at (0,1*\mainfact,1*\mainfact);
\coordinate [label=left:\text{$w_3=1,w_4=0)$}] (CC) at (0,0.93*\mainfact,1*\mainfact);
\coordinate [label=left:\text{$V_4=(w_1=0,w_2=0,w_3=0,w_4=1)$}] (D) at (0,0,0);

\draw[->] (0,0) -- (10,0,0) node[right] {$x$};
\draw[->] (0,0) -- (0,10,0) node[above] {$y$};
\draw[->] (0,0) -- (0,0,10) node[below left] {$z$};

\foreach \i in {A,B,C,D}
    \draw[dashed,very thick] (0,0)--(\i);
\draw[-,very thick, opacity=.5] (A)--(D)--(B)--cycle;
\draw[-,very thick, opacity=.5] (A) --(D)--(C)--cycle;
\draw[-,very thick, opacity=.5] (B)--(D)--(C)--cycle;
\draw[-,very thick] (A)--(B);
\draw[-,very thick] (A)--(C);
\draw[-,very thick] (A)--(D);

\coordinate [label=above:\text{}] (X1) at (1*\mainfact,1/2*\mainfact,1/2*\mainfact);

\draw[-,line width=0.6mm, opacity=.2, fill=brown!30] (C)--(X1)--(D)--cycle;
\draw[-,line width=0.6mm, dashed, draw=brown] (C)--(X1);
\draw[-,line width=0.6mm, draw=brown] (X1)--(D);

\coordinate [label=above:\text{}] (X3) at (7/8*\mainfact,7/8*\mainfact,0*\mainfact);

\draw[-,line width=0.6mm, opacity=.2, fill=violet!30] (B)--(X3)--(C)--cycle;
\draw[-,line width=0.6mm, draw=violet] (C)--(X3);
\draw[-,line width=0.6mm, draw=violet] (X3)--(B);

\coordinate [label=above:\text{}] (X4) at (2/3*\mainfact,1/3*\mainfact,1*\mainfact);

\draw[-,line width=0.6mm, opacity=.2, fill=orange!30] (A)--(X4)--(D)--cycle;
\draw[-,line width=0.6mm, dashed, draw=orange] (A)--(X4);
\draw[-,line width=0.6mm, dashed, draw=orange] (X4)--(D);

\coordinate [label=above:\text{}] (X6) at (0*\mainfact,1/4*\mainfact,1/4*\mainfact);

\draw[-,line width=0.6mm, opacity=.2, fill=magenta!30] (A)--(X6)--(B)--cycle;
\draw[-,line width=0.6mm, draw=magenta] (A)--(X6);
\draw[-,line width=0.6mm, dashed, draw=magenta] (X6)--(B);

\coordinate (E) at (0.5*\mainfact,0.375*\mainfact,0.375*\mainfact);
\coordinate (F) at (0.756756757*\mainfact,0.567567568*\mainfact,0.567567568*\mainfact);
\coordinate (G) at (0.913043478*\mainfact,0.47826087*\mainfact,0.47826087*\mainfact);
\coordinate (H) at (0.851851852*\mainfact,0.814814815*\mainfact,0.111111111*\mainfact);

The 1-2-3-4-1 cycle
\draw[-,very thick, fill=green!30, opacity=.5] (E)--(H)--(F)--cycle;
\draw[-,very thick, fill=green!30, opacity=.5] (E) --(H)--(G)--cycle;
\draw[-,very thick, fill=green!30, opacity=.5] (F)--(H)--(G)--cycle;

\tikzstyle{node1} = [circle,inner sep=2pt,draw=forestgreen,fill=forestgreen,thick];

\node [node1] (A1) at (8.5,2.5) {\footnotesize \color{white} 1};
\node [node1] (A2) at (10,2.5) {\footnotesize \color{white}  2};
\node [node1] (A3) at (10,1) {\footnotesize \color{white} 3};
\node [node1] (A4) at (8.5,1) {\footnotesize \color{white} 4};

  \draw [-latex,draw=forestgreen, very thick](A1) -- (A2);
  \draw [-latex,draw=forestgreen, very thick](A2) -- (A3);
  \draw [-latex,draw=forestgreen, very thick](A3) -- (A4);
  \draw [-latex,draw=forestgreen, very thick](A4) -- (A1);

\end{tikzpicture}

    \caption{The first tetrahedron of efficient weight vectors ($\mathcal{T}_1$), the associated Hamiltonian cycle, and the determining planes ($h_{12}, h_{14}, h_{23}$, and $h_{34}$) for the example.}
    \label{Cycle1Final}
\end{figure}

\begin{figure}[H]
    \centering

\begin{tikzpicture}[line join = round, line cap = round, scale=0.87]

\pgfmathsetmacro{\mainfact}{8};

\coordinate [label=above left:\text{$V_1=(w_1=1,w_2=0,w_3=0,w_4=0)$}] (A) at (1*\mainfact,1*\mainfact,0);
\coordinate [label=below:\text{$V_2=(w_1=0,w_2=1,w_3=0,w_4=0)$}] (B) at (1*\mainfact,0,1*\mainfact);
\coordinate [label=left:\text{$V_3=(w_1=0,w_2=0,$}] (C) at (0,1*\mainfact,1*\mainfact);
\coordinate [label=left:\text{$w_3=1,w_4=0)$}] (CC) at (0,0.93*\mainfact,1*\mainfact);
\coordinate [label=left:\text{$V_4=(w_1=0,w_2=0,w_3=0,w_4=1)$}] (D) at (0,0,0);

\draw[->] (0,0) -- (10,0,0) node[right] {$x$};
\draw[->] (0,0) -- (0,10,0) node[above] {$y$};
\draw[->] (0,0) -- (0,0,10) node[below left] {$z$};

\foreach \i in {A,B,C,D}
    \draw[dashed,very thick] (0,0)--(\i);
\draw[-,very thick, opacity=.5] (A)--(D)--(B)--cycle;
\draw[-,very thick, opacity=.5] (A) --(D)--(C)--cycle;
\draw[-,very thick, opacity=.5] (B)--(D)--(C)--cycle;
\draw[-,very thick] (A)--(B);
\draw[-,very thick] (A)--(C);
\draw[-,very thick] (A)--(D);

\coordinate [label=above:\text{}] (X2) at (5/6*\mainfact,1*\mainfact,1/6*\mainfact);

\draw[-,very thick, opacity=.2, fill=teal!30] (B)--(X2)--(D)--cycle;
\draw[-,line width=0.6mm, draw=teal] (D)--(X2);
\draw[-,line width=0.6mm, dashed, draw=teal] (X2)--(B);

\coordinate [label=above:\text{}] (X3) at (7/8*\mainfact,7/8*\mainfact,0*\mainfact);

\draw[-,line width=0.6mm, opacity=.2, fill=violet!30] (B)--(X3)--(C)--cycle;
\draw[-,line width=0.6mm, draw=violet] (C)--(X3);
\draw[-,line width=0.6mm, draw=violet] (X3)--(B);

\coordinate [label=above:\text{}] (X4) at (2/3*\mainfact,1/3*\mainfact,1*\mainfact);

\draw[-,line width=0.6mm, opacity=.2, fill=orange!30] (A)--(X4)--(D)--cycle;
\draw[-,line width=0.6mm, dashed, draw=orange] (A)--(X4);
\draw[-,line width=0.6mm, dashed, draw=orange] (X4)--(D);

\coordinate [label=above:\text{}] (X5) at (8/9*\mainfact,0*\mainfact,8/9*\mainfact);

\draw[-,line width=0.6mm, opacity=.2, fill=forestgreen!30] (A)--(X5)--(C)--cycle;
\draw[-,line width=0.6mm, draw=forestgreen] (A)--(X5);
\draw[-,line width=0.6mm, dashed, draw=forestgreen] (X5)--(C);

\coordinate (I) at (0.848484848*\mainfact,0.727272727*\mainfact,0.363636364*\mainfact);
\coordinate (J) at (0.803278689*\mainfact,0.68852459*\mainfact,0.344262295*\mainfact);
\coordinate (K) at (0.862068966*\mainfact,0.482758621*\mainfact,0.540229885*\mainfact);
\coordinate (L) at (0.75*\mainfact,0.55*\mainfact,0.6*\mainfact);

The 1-2-3-4-1 cycle
\draw[-,very thick, fill = blue!30, opacity =.5] (I)--(J)--(K)--cycle;
\draw[-,very thick, fill = blue!30, opacity =.5] (I) --(K)--(L)--cycle;
\draw[-,very thick, fill = blue!30, opacity =.5] (J)--(K)--(L)--cycle;

\tikzstyle{node2} = [circle,inner sep=2pt,draw=blue,fill= blue,thick];

\node [node2] (B1) at (8.5,2.5) {\footnotesize \color{white} 1};
\node [node2] (B2) at (8.5,1) {\footnotesize \color{white} 3};
\node [node2] (B3) at (10,1) {\footnotesize \color{white} 2};
\node [node2] (B4) at (10,2.5) {\footnotesize \color{white} 4};

  \draw [-latex,draw=blue, very thick](B2) -- (B1);
  \draw [-latex,draw=blue, very thick](B3) -- (B2);
  \draw [-latex,draw=blue, very thick](B4) -- (B3);
  \draw [-latex,draw=blue, very thick](B1) -- (B4);

\end{tikzpicture}

    \caption{The second tetrahedron of efficient weight vectors ($\mathcal{T}_2$), the associated Hamiltonian cycle, and the determining planes ($h_{13}, h_{14}, h_{23}$, and $h_{24}$) for the example.}
    \label{Cycle2Final}
\end{figure}

\begin{figure}[H]
    \centering

\begin{tikzpicture}[line join = round, line cap = round, scale=0.87]

\pgfmathsetmacro{\mainfact}{8};

\coordinate [label=above left:\text{$V_1=(w_1=1,w_2=0,w_3=0,w_4=0)$}] (A) at (1*\mainfact,1*\mainfact,0);
\coordinate [label=below:\text{$V_2=(w_1=0,w_2=1,w_3=0,w_4=0)$}] (B) at (1*\mainfact,0,1*\mainfact);
\coordinate [label=left:\text{$V_3=(w_1=0,w_2=0,$}] (C) at (0,1*\mainfact,1*\mainfact);
\coordinate [label=left:\text{$w_3=1,w_4=0)$}] (CC) at (0,0.93*\mainfact,1*\mainfact);
\coordinate [label=left:\text{$V_4=(w_1=0,w_2=0,w_3=0,w_4=1)$}] (D) at (0,0,0);

\draw[->] (0,0) -- (10,0,0) node[right] {$x$};
\draw[->] (0,0) -- (0,10,0) node[above] {$y$};
\draw[->] (0,0) -- (0,0,10) node[below left] {$z$};

\foreach \i in {A,B,C,D}
    \draw[dashed,very thick] (0,0)--(\i);
\draw[-,very thick, opacity=.5] (A)--(D)--(B)--cycle;
\draw[-,very thick, opacity=.5] (A) --(D)--(C)--cycle;
\draw[-,very thick, opacity=.5] (B)--(D)--(C)--cycle;
\draw[-,very thick] (A)--(B);
\draw[-,very thick] (A)--(C);
\draw[-,very thick] (A)--(D);

\coordinate [label=above:\text{}] (X1) at (1*\mainfact,1/2*\mainfact,1/2*\mainfact);

\draw[-,line width=0.6mm, opacity=.2, fill=brown!30] (C)--(X1)--(D)--cycle;
\draw[-,line width=0.6mm, dashed, draw=brown] (C)--(X1);
\draw[-,line width=0.6mm, draw=brown] (X1)--(D);

\coordinate [label=above:\text{}] (X2) at (5/6*\mainfact,1*\mainfact,1/6*\mainfact);

\draw[-,very thick, opacity=.2, fill=teal!30] (B)--(X2)--(D)--cycle;
\draw[-,line width=0.6mm, draw=teal] (D)--(X2);
\draw[-,line width=0.6mm, dashed, draw=teal] (X2)--(B);

\coordinate [label=above:\text{}] (X5) at (8/9*\mainfact,0*\mainfact,8/9*\mainfact);

\draw[-,line width=0.6mm, opacity=.2, fill=forestgreen!30] (A)--(X5)--(C)--cycle;
\draw[-,line width=0.6mm, draw=forestgreen] (A)--(X5);
\draw[-,line width=0.6mm, dashed, draw=forestgreen] (X5)--(C);

\coordinate [label=above:\text{}] (X6) at (0*\mainfact,1/4*\mainfact,1/4*\mainfact);

\draw[-,line width=0.6mm, opacity=.2, fill=magenta!30] (A)--(X6)--(B)--cycle;
\draw[-,line width=0.6mm, draw=magenta] (A)--(X6);
\draw[-,line width=0.6mm, dashed, draw=magenta] (X6)--(B);

\coordinate (M) at (0.923076923*\mainfact,0.480769231*\mainfact,0.480769231*\mainfact);
\coordinate (N) at (0.860215054*\mainfact,0.516129032*\mainfact,0.516129032*\mainfact);
\coordinate (O) at (0.714285714*\mainfact,0.428571429*\mainfact,0.428571429*\mainfact);
\coordinate (P) at (0.878787879*\mainfact,0.181818182*\mainfact,0.757575758*\mainfact);

The 1-2-3-4-1 cycle
\draw[-,very thick, fill = red!30, opacity =.5] (M)--(N)--(O)--cycle;
\draw[-,very thick, fill = red!30, opacity =.5] (M) --(O)--(P)--cycle;
\draw[-,very thick, fill = red!30, opacity =.5] (N)--(O)--(P)--cycle;

\tikzstyle{node3} = [circle,inner sep=2pt,draw=red,fill= red,thick];

\node [node3] (C1) at (8.5,2.5) {\footnotesize \color{white} 1};
\node [node3] (C2) at (8.5,1) {\footnotesize \color{white} 2};
\node [node3] (C3) at (10,1) {\footnotesize \color{white} 4};
\node [node3] (C4) at (10,2.5) {\footnotesize \color{white} 3};

  \draw [-latex,draw=red, very thick](C2) -- (C1);
  \draw [-latex,draw=red, very thick](C3) -- (C2);
  \draw [-latex,draw=red, very thick](C4) -- (C3);
  \draw [-latex,draw=red, very thick](C1) -- (C4);

\end{tikzpicture}

    \caption{The third tetrahedron of efficient weight vectors ($\mathcal{T}_3$), the associated Hamiltonian cycle, and the determining planes ($h_{12}, h_{13}, h_{24}$, and $h_{34}$) for the example.}
    \label{Cycle3Final}
\end{figure}

\begin{figure}[H]
    \centering

\begin{tikzpicture}[line join = round, line cap = round, scale=0.87]

\pgfmathsetmacro{\mainfact}{8};

\coordinate [label=above left:\text{$V_1=(w_1=1,w_2=0,w_3=0,w_4=0)$}] (A) at (1*\mainfact,1*\mainfact,0);
\coordinate [label=below:\text{$V_2=(w_1=0,w_2=1,w_3=0,w_4=0)$}] (B) at (1*\mainfact,0,1*\mainfact);
\coordinate [label=left:\text{$V_3=(w_1=0,w_2=0,$}] (C) at (0,1*\mainfact,1*\mainfact);
\coordinate [label=left:\text{$w_3=1,w_4=0)$}] (CC) at (0,0.93*\mainfact,1*\mainfact);
\coordinate [label=left:\text{$V_4=(w_1=0,w_2=0,w_3=0,w_4=1)$}] (D) at (0,0,0);

\draw[->] (0,0) -- (10,0,0) node[right] {$x$};
\draw[->] (0,0) -- (0,10,0) node[above] {$y$};
\draw[->] (0,0) -- (0,0,10) node[below left] {$z$};

\foreach \i in {A,B,C,D}
    \draw[dashed,very thick] (0,0)--(\i);
\draw[-,very thick, opacity=.5] (A)--(D)--(B)--cycle;
\draw[-,very thick, opacity=.5] (A) --(D)--(C)--cycle;
\draw[-,very thick, opacity=.5] (B)--(D)--(C)--cycle;
\draw[-,very thick] (A)--(B);
\draw[-,very thick] (A)--(C);
\draw[-,very thick] (A)--(D);

\coordinate (E) at (0.5*\mainfact,0.375*\mainfact,0.375*\mainfact);
\coordinate (F) at (0.756756757*\mainfact,0.567567568*\mainfact,0.567567568*\mainfact);
\coordinate (G) at (0.913043478*\mainfact,0.47826087*\mainfact,0.47826087*\mainfact);
\coordinate (H) at (0.851851852*\mainfact,0.814814815*\mainfact,0.111111111*\mainfact);

The 1-2-3-4-1 cycle
\draw[-,very thick, fill=green!30, opacity=.5] (E)--(H)--(F)--cycle;
\draw[-,very thick, fill=green!30, opacity=.5] (E) --(H)--(G)--cycle;
\draw[-,very thick, fill=green!30, opacity=.5] (F)--(H)--(G)--cycle;

\coordinate (I) at (0.848484848*\mainfact,0.727272727*\mainfact,0.363636364*\mainfact);
\coordinate (J) at (0.803278689*\mainfact,0.68852459*\mainfact,0.344262295*\mainfact);
\coordinate (K) at (0.862068966*\mainfact,0.482758621*\mainfact,0.540229885*\mainfact);
\coordinate (L) at (0.75*\mainfact,0.55*\mainfact,0.6*\mainfact);

The 1-2-3-4-1 cycle
\draw[-,very thick, fill = blue!30, opacity =.5] (I)--(J)--(K)--cycle;
\draw[-,very thick, fill = blue!30, opacity =.5] (I) --(K)--(L)--cycle;
\draw[-,very thick, fill = blue!30, opacity =.5] (J)--(K)--(L)--cycle;

\coordinate (M) at (0.923076923*\mainfact,0.480769231*\mainfact,0.480769231*\mainfact);
\coordinate (N) at (0.860215054*\mainfact,0.516129032*\mainfact,0.516129032*\mainfact);
\coordinate (O) at (0.714285714*\mainfact,0.428571429*\mainfact,0.428571429*\mainfact);
\coordinate (P) at (0.878787879*\mainfact,0.181818182*\mainfact,0.757575758*\mainfact);

The 1-2-3-4-1 cycle
\draw[-,very thick, fill = red!30, opacity =.5] (M)--(N)--(O)--cycle;
\draw[-,very thick, fill = red!30, opacity =.5] (M) --(O)--(P)--cycle;
\draw[-,very thick, fill = red!30, opacity =.5] (N)--(O)--(P)--cycle;

\tikzstyle{node1} = [circle,inner sep=2pt,draw=forestgreen,fill=forestgreen,thick];

\node [node1] (A1) at (8.5,7.5) {\footnotesize \color{white} 1};
\node [node1] (A2) at (10,7.5) {\footnotesize \color{white}  2};
\node [node1] (A3) at (10,6) {\footnotesize \color{white} 3};
\node [node1] (A4) at (8.5,6) {\footnotesize \color{white} 4};

  \draw [-latex,draw=forestgreen, very thick](A1) -- (A2);
  \draw [-latex,draw=forestgreen, very thick](A2) -- (A3);
  \draw [-latex,draw=forestgreen, very thick](A3) -- (A4);
  \draw [-latex,draw=forestgreen, very thick](A4) -- (A1);

  \tikzstyle{node2} = [circle,inner sep=2pt,draw=blue,fill= blue,thick];

\node [node2] (B1) at (8.5,5) {\footnotesize \color{white} 1};
\node [node2] (B2) at (8.5,3.5) {\footnotesize \color{white} 3};
\node [node2] (B3) at (10,3.5) {\footnotesize \color{white} 2};
\node [node2] (B4) at (10,5) {\footnotesize \color{white} 4};

  \draw [-latex,draw=blue, very thick](B2) -- (B1);
  \draw [-latex,draw=blue, very thick](B3) -- (B2);
  \draw [-latex,draw=blue, very thick](B4) -- (B3);
  \draw [-latex,draw=blue, very thick](B1) -- (B4);

\tikzstyle{node3} = [circle,inner sep=2pt,draw=red,fill= red,thick];

\node [node3] (C1) at (8.5,2.5) {\footnotesize \color{white} 1};
\node [node3] (C2) at (8.5,1) {\footnotesize \color{white} 2};
\node [node3] (C3) at (10,1) {\footnotesize \color{white} 4};
\node [node3] (C4) at (10,2.5) {\footnotesize \color{white} 3};

  \draw [-latex,draw=red, very thick](C2) -- (C1);
  \draw [-latex,draw=red, very thick](C3) -- (C2);
  \draw [-latex,draw=red, very thick](C4) -- (C3);
  \draw [-latex,draw=red, very thick](C1) -- (C4);

\end{tikzpicture}

    \caption{The set of efficient weight vectors for the example (the union of $\mathcal{T}_1, \mathcal{T}_2$ and $\mathcal{T}_3$). An interactive version is available at https://www.geogebra.org/m/hnk47sz5.}
    \label{OnlyOptFinal}
\end{figure}

\end{example}

\section{Conclusion and further research}
\label{sec:4}

In this paper the first formal proof on the geometry of Pareto efficient weight vectors is provided. It has been shown that the set of efficient weight vectors for a $4\times4$ pairwise comparison matrix $\mathbf{A}$ can be represented as the union of three tetrahedra. Each tetrahedron is determined by four weight vectors calculated from four incomplete submatrices of $\mathbf{A}$ that can be represented by path spanning tree graphs related to a given $4$-cycle. It is also proved that with appropriate rearrangements, the orientations of $4$-cycles in the Blanquero-Carrozisa-Conde graphs (which determine efficiency) are also fixed.

From a practical point of view it is shown that a weight vector is efficient for a $4\times4$ pairwise comparison matrix $\mathbf{A}$ if and only if it is a convex combination of four weight vectors calculated from path spanning trees that determine exactly one of the $(1,2,3,4,1)$, $(1,4,2,3,1)$ and $(1,3,4,2,1)$ cycles. All possible special cases have been illustrated regarding the geometry of the set of efficient weight vectors in the main text of the paper and in Appendix~\ref{append:B}.

The geometry of efficient weight vectors is relevant for larger pairwise comparison matrices as well, although visualization has already achieved its limit at four alternatives ($n=4$). The pattern shown above suggests that each of the $\frac{1}{2}(n-1)!$ Hamiltonian cycles of the BCC graph corresponds to a simplex, determined by $n$ spanning trees' weight vectors. These simplices form building blocks of the set of efficient weight vectors. 
Within each simplex, the weight vectors associated with the paths of length  $n-1$  form a local basis. 
Every efficient weight vector can be written as a convex combination of $n$ spanning trees' weight vectors, belonging to the same Hamiltonian cycle. Note that all the other spanning trees' (different from a path of length $n-1$) weight vectors can also be written in that way. Generalizations of Propositions~\ref{negykor} and \ref{final} are also necessary in future research.

\section*{Acknowledgements}
The research was supported by the National Research, Development and Innovation Office under Grants FK 145838 and TKP2021-NKTA-01 NRDIO.

\bibliographystyle{apalike} 
\bibliography{main}

\begin{thebibliography}{}

\bibitem[\'Abele{-}Nagy and Bozóki, 2016]{Abele2016}
\'Abele{-}Nagy, K. and Bozóki, S. (2016).
\newblock Efficiency analysis of simple perturbed pairwise comparison matrices.
\newblock {\em Fundamenta Informaticae}, 144(3-4):279--289.
\newblock \url{https://doi.org/10.3233/FI-2016-1335}.

\bibitem[\'Abele{-}Nagy et~al., 2018]{Abele2018}
\'Abele{-}Nagy, K., Bozóki, S., and Rebák, {\"O}. (2018).
\newblock Efficiency analysis of double perturbed pairwise comparison matrices.
\newblock {\em Journal of the Operational Research Society}, 69(5):707--713.
\newblock \url{https://doi.org/10.1080/01605682.2017.1409408}.

\bibitem[Bajwa et~al., 2008]{BajwaChooWedley2008}
Bajwa, G., Choo, E.~U., and Wedley, W.~C. (2008).
\newblock Effectiveness analysis of deriving priority vectors from reciprocal pairwise comparison matrices.
\newblock {\em Asia-Pacific Journal of Operational Research}, 25(3):279--299.
\newblock \url{https://doi.org/10.1142/S0217595908001754}.

\bibitem[Blanquero et~al., 2006]{Blanquero}
Blanquero, R., Carrizosa, E., and Conde, E. (2006).
\newblock Inferring efficient weights from pairwise comparison matrices.
\newblock {\em Mathematical Methods of Operations Research}, 64(2):271--284.
\newblock \url{https://doi.org/10.1007/s00186-006-0077-1}.

\bibitem[Bozóki, 2014]{Bozoki2014}
Bozóki, S. (2014).
\newblock Inefficient weights from pairwise comparison matrices with arbitrarily small inconsistency.
\newblock {\em Optimization}, 63(12):1893--1901.
\newblock \url{https://doi.org/10.1080/02331934.2014.903399}.

\bibitem[Bozóki and Fülöp, 2018]{BozokiFulop2018}
Bozóki, S. and Fülöp, J. (2018).
\newblock Efficient weight vectors from pairwise comparison matrices.
\newblock {\em European Journal of Operational Research}, 264(2):419--427.
\newblock \url{https://doi.org/10.1016/j.ejor.2017.06.033}.

\bibitem[Bozóki et~al., 2010]{Bozoki2010}
Bozóki, S., Fülöp, J., and Rónyai, L. (2010).
\newblock On optimal completion of incomplete pairwise comparison matrices.
\newblock {\em Mathematical and Computer Modelling}, 52(1):318--333.
\newblock \url{https://doi.org/10.1016/j.mcm.2010.02.047}.

\bibitem[Brunelli, 2018]{Brunelli2018}
Brunelli, M. (2018).
\newblock A survey of inconsistency indices for pairwise comparisons.
\newblock {\em International Journal of General Systems}, 47(8):751--771.
\newblock \url{https://doi.org/10.1080/03081079.2018.1523156}.

\bibitem[Brunelli and Fedrizzi, 2024]{BrunelliFedrizzi2024}
Brunelli, M. and Fedrizzi, M. (2024).
\newblock Inconsistency indices for pairwise comparisons and the {P}areto dominance principle.
\newblock {\em European Journal of Operational Research}, 312(1):273--282.
\newblock \url{https://doi.org/10.1016/j.ejor.2023.06.033}.

\bibitem[Camion, 1959]{Camion1959}
Camion, P. (1959).
\newblock Chemins et circuits hamiltoniens des graphes complets.
\newblock {\em Comptes Rendus de l'Académie des Sciences de Paris}, 249:2151--2152.
\newblock In French. \url{https://gallica.bnf.fr/ark:/12148/bpt6k731d/f1025.item}.

\bibitem[Choo and Wedley, 2004]{ChooWedley2004}
Choo, E. and Wedley, W. (2004).
\newblock A common framework for deriving preference values from pairwise comparison matrices.
\newblock {\em Computers \& Operations Research}, 31(6):893--908.
\newblock \url{https://doi.org/10.1016/S0305-0548(03)00042-X}.

\bibitem[Conde and {de la Paz Rivera Pérez}, 2010]{Conde2010}
Conde, E. and {de la Paz Rivera Pérez}, M. (2010).
\newblock A linear optimization problem to derive relative weights using an interval judgement matrix.
\newblock {\em European Journal of Operational Research}, 201(2):537--544.
\newblock \url{https://doi.org/10.1016/j.ejor.2009.03.029}.

\bibitem[Crawford and Williams, 1985]{Crawford1985}
Crawford, G. and Williams, C. (1985).
\newblock A note on the analysis of subjective judgment matrices.
\newblock {\em Journal of Mathematical Psychology}, 29(4):387--405.
\newblock \url{https://doi.org/10.1016/0022-2496(85)90002-1}.

\bibitem[Csató, 2021]{Csato2021}
Csató, L. (2021).
\newblock {\em Tournament Design: How Operations Research Can Improve Sports Rules}.
\newblock Palgrave Pivots in Sports Economics, Palgrave Macmillan.
\newblock \url{https://doi.org/10.1007/978-3-030-59844-0}.

\bibitem[{da Cruz} et~al., 2021]{DaCruz2021}
{da Cruz}, H.~F., Fernandes, R., and Furtado, S. (2021).
\newblock Efficient vectors for simple perturbed consistent matrices.
\newblock {\em International Journal of Approximate Reasoning}, 139:54--68.
\newblock \url{https://doi.org/10.1016/j.ijar.2021.09.007}.

\bibitem[Dijkstra, 2013]{Dijkstra2013}
Dijkstra, T.~K. (2013).
\newblock On the extraction of weights from pairwise comparison matrices.
\newblock {\em Central European Journal of Operations Research}, 21:103--123.
\newblock \url{https://doi.org/10.1007/s10100-011-0212-9}.

\bibitem[Duleba and Moslem, 2019]{DulebaMoslem2019}
Duleba, {\relax Sz}. and Moslem, S. (2019).
\newblock Examining {P}areto optimality in analytic hierarchy process on real {D}ata: An application in public transport service development.
\newblock {\em Expert Systems with Applications}, 116:21--30.
\newblock \url{https://doi.org/10.1016/j.eswa.2018.08.049}.

\bibitem[Ehrgott, 2000]{Ehrgott2000}
Ehrgott, M. (2000).
\newblock {\em Multicriteria Optimization}.
\newblock volume 491 of Lecture Notes in Economics and Mathematical Systems. Springer Verlag, Berlin.

\bibitem[Fernandes and Furtado, 2022]{Fernandes2022}
Fernandes, R. and Furtado, S. (2022).
\newblock Efficiency of the principal eigenvector of some triple perturbed consistent matrices.
\newblock {\em European Journal of Operational Research}, 298(3):1007--1015.
\newblock \url{https://doi.org/10.1016/j.ejor.2021.08.012}.

\bibitem[Furtado, 2023]{Furtado2023}
Furtado, S. (2023).
\newblock Efficient vectors for double perturbed consistent matrices.
\newblock {\em Optimization}, 72(11):2679--2701.
\newblock \url{ https://doi.org/10.1080/02331934.2022.2070067}.

\bibitem[Furtado and Johnson, 2024a]{FurtadoJohnson2024a}
Furtado, S. and Johnson, C. (2024a).
\newblock Efficient vectors for block perturbed consistent matrices.
\newblock {\em SIAM Journal on Matrix Analysis and Applications}, 45(1):601--618.
\newblock \url{ https://doi.org/10.1137/23M1580310}.

\bibitem[Furtado and Johnson, 2024b]{FurtadoJohnson2024b}
Furtado, S. and Johnson, C.~R. (2024b).
\newblock Efficiency analysis for the {P}erron vector of a reciprocal matrix.
\newblock {\em Applied Mathematics and Computation}, 480:128913.
\newblock \url{https://doi.org/10.1016/j.amc.2024.128913}.

\bibitem[Furtado and Johnson, 2024c]{FurtadoJohnson2024c}
Furtado, S. and Johnson, C.~R. (2024c).
\newblock Efficiency of any weighted geometric mean of the columns of a reciprocal matrix.
\newblock {\em Linear Algebra and its Applications}, 680:83--92.
\newblock \url{https://doi.org/10.1016/j.laa.2023.10.001}.

\bibitem[Furtado and Johnson, 2024d]{FurtadoJohnson2024d}
Furtado, S. and Johnson, C.~R. (2024d).
\newblock Efficiency of the convex hull of the columns of certain triple perturbed consistent matrices.
\newblock {\em Advances in Operator Theory}, 9(4):89.
\newblock \url{https://doi.org/10.1007/s43036-024-00384-z}.

\bibitem[Furtado and Johnson, 2024e]{FurtadoJohnson2024e}
Furtado, S. and Johnson, C.~R. (2024e).
\newblock Efficient vectors in priority setting methodology.
\newblock {\em Annals of Operations Research}, 332:743--764.
\newblock \url{ https://doi.org/10.1007/s10479-023-05771-y}.

\bibitem[Furtado and Johnson, 2024f]{FurtadoJohnson2024f}
Furtado, S. and Johnson, C.~R. (2024f).
\newblock Pairwise comparison matrices with uniformly ordered efficient vectors.
\newblock {\em International Journal of Approximate Reasoning}, 173:109265.
\newblock \url{https://doi.org/10.1016/j.ijar.2024.109265}.

\bibitem[Furtado and Johnson, 2025]{FurtadoJohnson-completeset-2025}
Furtado, S. and Johnson, C.~R. (2025).
\newblock The complete set of efficient vectors for a pairwise comparison matrix.
\newblock {\em Computational and Applied Mathematics}, 44:100.
\newblock \url{https://doi.org/10.1007/s40314-024-03062-1}.

\bibitem[Gass, 1998]{Gass1998}
Gass, S. (1998).
\newblock Tournaments, transitivity and pairwise comparison matrices.
\newblock {\em Journal of the Operational Research Society}, 49(6):616--624.
\newblock \url{https://doi.org/10.1057/palgrave.jors.2600572}.

\bibitem[Golany and Kress, 1993]{GolanyKress1993}
Golany, B. and Kress, M. (1993).
\newblock A multicriteria evaluation of methods for obtaining weights from ratio-scale matrices.
\newblock {\em European Journal of Operational Research}, 69(2):210--220.
\newblock \url{https://doi.org/10.1016/0377-2217(93)90165-J}.

\bibitem[Harker, 1987]{Harker1987}
Harker, P.~T. (1987).
\newblock Incomplete pairwise comparisons in the analytic hierarchy process.
\newblock {\em Mathematical Modelling}, 9(11):837--848.
\newblock \url{https://doi.org/10.1016/0270-0255(87)90503-3}.

\bibitem[Jensen, 1984]{Jensen1984}
Jensen, R.~E. (1984).
\newblock An alternative scaling method for priorities in hierarchical structures.
\newblock {\em Journal of Mathematical Psychology}, 28(3):317--332.
\newblock \url{https://doi.org/10.1016/0022-2496(84)90003-8}.

\bibitem[Lundy et~al., 2017]{Lundy2017}
Lundy, M., Siraj, S., and Greco, S. (2017).
\newblock The mathematical equivalence of the “spanning tree” and row geometric mean preference vectors and its implications for preference analysis.
\newblock {\em European Journal of Operational Research}, 257(1):197--208.
\newblock \url{https://doi.org/10.1016/j.ejor.2016.07.042}.

\bibitem[Mazurek and Kułakowski, 2022]{Mazurek2022}
Mazurek, J. and Kułakowski, K. (2022).
\newblock On the derivation of weights from incomplete pairwise comparisons matrices via spanning trees with crisp and fuzzy confidence levels.
\newblock {\em International Journal of Approximate Reasoning}, 150:242--257.
\newblock \url{https://doi.org/10.1016/j.ijar.2022.08.014}.

\bibitem[Saaty, 1977]{Saaty1977}
Saaty, T.~L. (1977).
\newblock A scaling method for priorities in hierarchical structures.
\newblock {\em Journal of Mathematical Psychology}, 15(3):234--281.
\newblock \url{https://doi.org/10.1016/0022-2496(77)90033-5}.

\bibitem[Saaty, 1980]{Saaty1980}
Saaty, T.~L. (1980).
\newblock {\em The Analytic Hierarchy Process}.
\newblock McGraw-Hill, New York.

\bibitem[Steuer, 1986]{Steuer1986}
Steuer, R.~E. (1986).
\newblock {\em Multiple Criteria Optimization: Theory, Computation, and Application}.
\newblock Wiley Series in Probability and Mathematical Statistics. Wiley.

\bibitem[Szádoczki and Bozóki, 2024]{Szadoczki2024}
Szádoczki, {\relax Zs}. and Bozóki, S. (2024).
\newblock Geometric interpretation of efficient weight vectors.
\newblock {\em Knowledge-Based Systems}, 303:112403.
\newblock \url{https://doi.org/10.1016/j.knosys.2024.112403}.

\bibitem[Tekile et~al., 2023]{Tekile2023}
Tekile, H.~A., Brunelli, M., and Fedrizzi, M. (2023).
\newblock A numerical comparative study of completion methods for pairwise comparison matrices.
\newblock {\em Operations Research Perspectives}, 10:100272.
\newblock \url{https://doi.org/10.1016/j.orp.2023.100272}.

\bibitem[Thurstone, 1927]{Thurstone1927}
Thurstone, L. (1927).
\newblock A law of comparative judgment.
\newblock {\em Psychological Review}, 34(4):273–286.
\newblock \url{https://doi.org/10.1037/h0070288}.

\bibitem[Triantaphyllou, 2000]{Triantaphyllou2000}
Triantaphyllou, E. (2000).
\newblock Multi-criteria decision making methods.
\newblock In {\em Multi-criteria Decision Making Methods: A Comparative Study}. Applied Optimization, vol 44. Springer, Boston, MA.
\newblock \url{https://doi.org/10.1007/978-1-4757-3157-6_2}.

\bibitem[Tsyganok, 2010]{Tsyganok2010}
Tsyganok, V. (2010).
\newblock Investigation of the aggregation effectiveness of expert estimates obtained by the pairwise comparison method.
\newblock {\em Mathematical and Computer Modelling}, 52(3):538--544.
\newblock \url{https://doi.org/10.1016/j.mcm.2010.03.052}.

\end{thebibliography}
\addcontentsline{toc}{section}{References}

\newpage

\newcounter{appendproposition}
\setcounter{appendproposition}{1}
\newcounter{appendexample}
\setcounter{appendexample}{1}
\newcounter{appendfigure}
\setcounter{appendfigure}{1}
\newcounter{appendremark}
\setcounter{appendremark}{1}
\newcounter{append}
\renewcommand{\theappend}{A}
\section*{Appendix A}
\addcontentsline{toc}{section}{Appendix A}
\refstepcounter{append}
\label{append:A}

If the triads of a $4 \times 4$ PCM determine its properties, then Proposition~\ref{atrendezes} can be used to categorize these matrices.

\renewcommand{\theproposition}{A\arabic{appendproposition}}
\begin{proposition} \label{atrendezes}
Any $4\times 4$ pairwise comparison matrix 
\[ \mathbf{A} =
\begin{pmatrix}
1 & a_{12} & a_{13} & a_{14} \\
1/a_{12} & 1 & a_{23} & a_{24} \\
1/a_{13} & 1/a_{23} & 1 & a_{34} \\
1/a_{14} & 1/a_{24} & 1/a_{34} & 1
\end{pmatrix},
\]
for which Assumption~\ref{assumption:1} holds, can be transformed with appropriate reindexing of alternatives to the following form:
\[ \mathbf{A}' =
\begin{pmatrix}
1 & a_{12}'& a_{13}' & a_{14}' \\
1/a_{12}' & 1 & a_{23}' & a_{24}' \\
1/a_{13}' & 1/a_{23}' & 1 & a_{34}' \\
1/a_{14}' & 1/a_{24}' & 1/a_{34}' & 1
\end{pmatrix},
\]
where $\mathbf{A}'$ falls into one of the following cases:

\bigskip

\noindent Case 1: \quad $a_{12}'a_{23}'<a_{13}'$, \quad $a_{23}'a_{34}'<a_{24}'$, \quad  $a_{13}'a_{34}'<a_{14}'$, \quad  $a_{12}'a_{24}'<a_{14}'$. \\
Case 2: \quad $a_{12}'a_{23}'<a_{13}'$, \quad $a_{23}'a_{34}'<a_{24}'$, \quad  $a_{13}'a_{34}'<a_{14}'$, \quad  $a_{12}'a_{24}'>a_{14}'$.
\end{proposition}
\begin{proof}
Let us consider $\mathbf{A}$ and let us swap the index of the first two alternatives. This will result in the following pairwise comparison matrix:
\[ \mathbf{A_{21}} =
\begin{pmatrix}
1 & 1/a_{12} & a_{23} & a_{24} \\
a_{12} & 1 & a_{13} & a_{14} \\
1/a_{23} & 1/a_{13} & 1 & a_{34} \\
1/a_{24} & 1/a_{14} & 1/a_{34} & 1
\end{pmatrix}
\]
Let us assume that in $\mathbf{A}$ 
\begin{itemize}
\item for triad (1,2,3) $a_{12}a_{23}>a_{13}$, 
\item for triad (2,3,4) $a_{23}a_{34}>a_{24}$, 
\item for triad (1,3,4) $a_{13}a_{34}>a_{14}$,
\item for triad (1,2,4) $a_{12}a_{24}>a_{14}$.
\end{itemize}

Then, in $\mathbf{A_{21}}$
\begin{itemize}
\item for triad (1,2,3) $(1/a_{12}) \cdot a_{13} < a_{23}$, 
\item for triad (2,3,4) $a_{13}a_{34}>a_{14}$,  
\item for triad (1,3,4) $a_{23}a_{34}>a_{24}$,
\item for triad (1,2,4) $(1/a_{12}) \cdot a_{14} < a_{24}$.
\end{itemize}

This transformation reversed the direction of the relation in those triads where both the reindexed alternatives were involved, and left them in place where only one of them was involved. With the same reasoning, a further swap of the indices of alternatives 3 and 4 will result in the form in Case 1. This is true if any pair of triads have their relation in a $>$ direction, swapping their indices will result in the relations corresponding to those triads to become $<$. Thus, if exactly zero, two or four of the triads have their relation in a $>$ direction, performing zero, one, or two index swappings of the appropriate alternatives will result in the form of Case 1.

On the other hand, if one or three triads have their relation in a $>$ direction, performing an appropriate number of index swaps will result in the form of Case 2.
\end{proof}

This means that if the triads of a $4\times 4$ matrix determine its properties, then it is enough to examine these two cases, and show the geometric nature of the set of efficient weight vectors, and the directions of the arcs in the strongly connected BCC graphs. 

However, Example~\ref{ellenpelda} shows that this is not the case.

\renewcommand{\theexample}{A\arabic{appendexample}}
\begin{example} \label{ellenpelda} Let us consider the following PCMs:
\[\underline{\mathbf{A}} =
\begin{pmatrix}
1       & 1 & 2   & 4       \\
1 &     1     &  1 & 3        \\
\frac{1}{2} &    1      &  1   & 1 \\
\frac{1}{4} &     \frac{1}{3}     &      1  & 1
\end{pmatrix},
\quad \bar{\mathbf{A}} =
\begin{pmatrix}
1       & 1 & 2   & \color{red} 8       \\
1 &     1     &  1 & 3        \\
\frac{1}{2} &    1      &  1   & 1 \\
\color{red} \frac{1}{8} &     \frac{1}{3}     &      1  & 1
\end{pmatrix},
\quad \mathbf{A} =
\begin{pmatrix}
1       & 1 & 2   & \color{red} 6       \\
1 &     1     &  1 & 3        \\
\frac{1}{2} &    1      &  1   & 1 \\
\color{red} \frac{1}{6} &     \frac{1}{3}     &      1  & 1
\end{pmatrix}.
\]

All of these matrices fall into Case 1 from Proposition~\ref{atrendezes}, as the relations hold even with the smallest $a_{14}=4$ value: 
\bigskip

\noindent $a_{12}'a_{23}'<a_{13}'$, \quad $a_{23}'a_{34}'<a_{24}'$, \quad  $a_{13}'a_{34}'<a_{14}'$, \quad  $a_{12}'a_{24}'<a_{14}'$. \\
 \qquad \ $1\cdot 1<2$, \qquad \ \ $1\cdot 1<3$, \qquad \ \ $2\cdot1<4$, \qquad  \ $1\cdot 3<4$.

 Thus, if the triads characterize the orientation of the Hamiltonian cycles associated with $\mathcal{T}_1,\mathcal{T}_2,$ and $\mathcal{T}_3$, then those should be the same for all of these matrices. However, this is not the case, Figure~\ref{Examplecycles} shows the cycles associated with each tetrahedron for $\underline{\mathbf{A}}$ and $\bar{\mathbf{A}}$.

 \renewcommand{\thefigure}{A\arabic{appendfigure}}
\begin{figure}[H]
    \centering

\begin{tikzpicture}[line join = round, line cap = round]

\tikzstyle{node1} = [circle,inner sep=2pt,draw=black,fill=black,thick];

\node [node1] (A1) at (4.5,-1) {\footnotesize \color{white} 1};
\node [node1] (A2) at (6,-1) {\footnotesize \color{white}  2};
\node [node1] (A3) at (6,-2.5) {\footnotesize \color{white} 3};
\node [node1] (A4) at (4.5,-2.5) {\footnotesize \color{white} 4};

  \draw [-latex,draw=black, very thick](A1) -- (A2);
  \draw [-latex,draw=black, very thick](A2) -- (A3);
  \draw [-latex,draw=black, very thick](A3) -- (A4);
  \draw [-latex,draw=black, very thick](A4) -- (A1);

\tikzstyle{node2} = [circle,inner sep=2pt,draw=black,fill= black,thick];

\node [node2] (B1) at (8.5,-1) {\footnotesize \color{white} 1};
\node [node2] (B2) at (8.5,-2.5) {\footnotesize \color{white} 3};
\node [node2] (B3) at (10,-2.5) {\footnotesize \color{white} 2};
\node [node2] (B4) at (10,-1) {\footnotesize
 \color{white} 4};

  \draw [-latex,draw=black, very thick](B2) -- (B1);
  \draw [-latex,draw=black, very thick](B3) -- (B2);
  \draw [-latex,draw=black, very thick](B4) -- (B3);
  \draw [-latex,draw=black, very thick](B1) -- (B4);

\tikzstyle{node3} = [circle,inner sep=2pt,draw=black,fill= black,thick];

\node [node3] (C1) at (12.5,-1) {\footnotesize \color{white} 1};
\node [node3] (C2) at (12.5,-2.5) {\footnotesize \color{white} 2};
\node [node3] (C3) at (14,-2.5) {\footnotesize \color{white} 4};
\node [node3] (C4) at (14,-1) {\footnotesize \color{white} 3};

  \draw [-latex,draw=black, very thick](C2) -- (C1);
  \draw [-latex,draw=black, very thick](C3) -- (C2);
  \draw [-latex,draw=black, very thick](C4) -- (C3);
  \draw [-latex,draw=black, very thick](C1) -- (C4);


\node [node1] (AA1) at (4.5,-3.5) {\footnotesize \color{white} 1};
\node [node1] (AA2) at (6,-3.5) {\footnotesize \color{white}  2};
\node [node1] (AA3) at (6,-5) {\footnotesize \color{white} 3};
\node [node1] (AA4) at (4.5,-5) {\footnotesize \color{white} 4};

  \draw [-latex,draw=black, very thick](AA1) -- (AA2);
  \draw [-latex,draw=black, very thick](AA2) -- (AA3);
  \draw [-latex,draw=black, very thick](AA3) -- (AA4);
  \draw [-latex,draw=black, very thick](AA4) -- (AA1);

\tikzstyle{node2} = [circle,inner sep=2pt,draw=black,fill= black,thick];

\node [node2] (BB1) at (8.5,-3.5) {\footnotesize \color{white} 1};
\node [node2] (BB2) at (8.5,-5) {\footnotesize \color{white} 3};
\node [node2] (BB3) at (10,-5) {\footnotesize \color{white} 2};
\node [node2] (BB4) at (10,-3.5) {\footnotesize
 \color{white} 4};

  \draw [-latex,draw=black, very thick](BB1) -- (BB2);
  \draw [-latex,draw=black, very thick](BB2) -- (BB3);
  \draw [-latex,draw=black, very thick](BB3) -- (BB4);
  \draw [-latex,draw=black, very thick](BB4) -- (BB1);

\tikzstyle{node3} = [circle,inner sep=2pt,draw=black,fill= black,thick];

\node [node3] (CC1) at (12.5,-3.5) {\footnotesize \color{white} 1};
\node [node3] (CC2) at (12.5,-5) {\footnotesize \color{white} 2};
\node [node3] (CC3) at (14,-5) {\footnotesize \color{white} 4};
\node [node3] (CC4) at (14,-3.5) {\footnotesize \color{white} 3};

  \draw [-latex,draw=black, very thick](CC2) -- (CC1);
  \draw [-latex,draw=black, very thick](CC3) -- (CC2);
  \draw [-latex,draw=black, very thick](CC4) -- (CC3);
  \draw [-latex,draw=black, very thick](CC1) -- (CC4);


\node (T1) at (5.25,0) {\footnotesize $\mathcal{T}_1$};
\node (T2) at (9.25,0) {\footnotesize $\mathcal{T}_2$};
\node (T3) at (13.25,0) {\footnotesize $\mathcal{T}_3$};


\node (Au) at (3.5,-1.75) {\footnotesize $\underline{\mathbf{A}}$};
\node (Ab) at (3.5,-4.25) {\footnotesize $\bar{\mathbf{A}}$};

\end{tikzpicture}

    \caption{The labeled Hamiltonian cycles of BCC graphs associated with tetrahedra $\mathcal{T}_1,\mathcal{T}_2,$ and $\mathcal{T}_3$ for $\underline{\mathbf{A}}$ and $\bar{\mathbf{A}}$. Each tetrahedron is related to the cycles below them, and each row belongs to the matrix at the beginning of the row.}
    \label{Examplecycles}
\end{figure}

 Note that the orientation of the Hamiltonian cycle associated with $\mathcal{T}_2$ is different for $\underline{\mathbf{A}}$ and $\bar{\mathbf{A}}$. The reason can be shown by $\mathbf{A}$, because that is the turning point between the two orientations. As all of these matrices fall into Case 1, that means all of the triads are inconsistent. However, the 4-cycle of $(1,3,2,4,1)$ is consistent for matrix $\mathbf{A}$, as $a_{13}a_{32}a_{24}a_{41}=1$ ($2\cdot1\cdot3\cdot1/6=1$). This also means that $\mathbf{A}$ is a double perturbed matrix, changing two elements (i.e., let $a_{12}=2$ and $a_{34}=3$) makes it consistent. As we differ from the double perturbed case (the consistent 4-cycle) in one ($\underline{\mathbf{A}}$) or the other ($\bar{\mathbf{A}}$) direction, the orientation of the Hamiltonian cycle associated with $\mathcal{T}_2$ that determines a subset of the efficient weight vectors changes as well.
\end{example}

\renewcommand{\theremark}{A\arabic{appendremark}}
\begin{remark} \label{prop1case2}
    All of the matrices in Case 2 provide the same orientations of the Hamiltonian cycles associated with $\mathcal{T}_1,\mathcal{T}_2,$ and $\mathcal{T}_3$ as $\underline{\mathbf{A}}$.
\end{remark}

\newpage

\newcounter{appendbproposition}
\setcounter{appendbproposition}{1}
\newcounter{appendbexample}
\setcounter{appendbexample}{1}
\newcounter{appendbfigure}
\setcounter{appendbfigure}{1}
\newcounter{appendbremark}
\setcounter{appendbremark}{1}
\newcounter{appendb}
\renewcommand{\theappendb}{B}
\section*{Appendix B}
\addcontentsline{toc}{section}{Appendix B}
\refstepcounter{appendb}
\label{append:B}

For the sake of brevity, whenever a given element and/or its reciprocal should be mentioned, we will only refer to the element.

\renewcommand{\theremark}{B\arabic{appendbremark}}

\begin{remark}
    In any $4\times4$ pairwise comparison matrix $\mathbf{A}$ there are the following four triads (undirected $3$-cycles), where $\square$ can be any of $<$, $>$, or $=$.

\begin{itemize}
    \item $(1,2,3): a_{12}a_{23}\square a_{13} \iff a_{12}a_{23}a_{31}\square 1,$
    \item $(1,2,4): a_{12}a_{24}\square a_{14} \iff a_{12}a_{24}a_{41}\square 1,$
    \item $(1,3,4): a_{13}a_{34}\square a_{14} \iff a_{13}a_{34}a_{41}\square 1,$
    \item $(2,3,4): a_{23}a_{34}\square a_{24} \iff a_{23}a_{34}a_{42}\square 1$.
\end{itemize}

There are altogether the following three (undirected) $4$-cycles with the same notations.

\begin{itemize}
    \item $(1,2,3,4): a_{12}a_{23}a_{34}a_{41}\square 1,$
    \item $(1,4,2,3): a_{14}a_{42}a_{23}a_{31}\square 1,$
    \item $(1,3,4,2): a_{13}a_{34}a_{42}a_{21}\square 1$.
\end{itemize}

Note the following properties and see Figure~\ref{SpecCasesCycles} for illustration:

\begin{enumerate}[label=(\roman*)]
    \item \label{prop:1} each element of the PCM appears in exactly two triads;
    \item \label{prop:2} each element of the PCM appears in exactly two $4$-cycles;
    \item \label{prop:3} for a given triad there is exactly one common element with each of the other triads;
    \item \label{prop:4} for a given triad there are exactly two common elements with each of the $4$-cycles,
    \item \label{prop:5} for a given $4$-cycle there are exactly two common elements with each of the other $4$-cycles,
    \item \label{prop:6} the common elements of any two $4$-cycles do not appear in the same triad,
    \item \label{prop:7} the different elements of any two $4$-cycles do not appear in the same triad,
    \item \label{prop:8} however two triads are chosen, those determine exactly one $4$-cycle, there is only one element that is not included in their union, and that element is included in both of the other two $4$-cycles,
    \item \label{prop:9} however two elements are chosen that appear in a common triad (for any two edges with a common starting vertex in Figure~\ref{SpecCasesCycles}) at least one of them will be included in all possible $4$-cycles.
\end{enumerate}
\end{remark}

\renewcommand{\thefigure}{B\arabic{appendbfigure}}
\begin{figure}[H]
    \centering

\begin{tikzpicture}[line join = round, line cap = round]

\tikzstyle{nodeS} = [circle,inner sep=2pt,draw=black,ultra thick];

\node [nodeS] (A1) at (2,-1) {\footnotesize 1};
\node [nodeS] (A2) at (5,-1) {\footnotesize 2};
\node [nodeS] (A3) at (5,-4) {\footnotesize 3};
\node [nodeS] (A4) at (2,-4) {\footnotesize 4};

  \draw [draw=black, very thick](A1) -- (A2);
  \draw [draw=black, very thick](A1) -- (A3);
  \draw [draw=black, very thick](A2) -- (A3);
  \draw [draw=black, very thick](A2) -- (A4);
  \draw [draw=black, very thick](A3) -- (A4);
  \draw [draw=black, very thick](A4) -- (A1);

  \draw (A1) -- (A2) node [midway, above=2pt] {$a_{12}$};
  \draw (A2) -- (A3) node [midway, right=2pt] {$a_{23}$};
  \draw (A3) -- (A4) node [midway, below=2pt] {$a_{34}$};
  \draw (A4) -- (A1) node [midway, left=2pt] {$a_{14}$};
  \draw (A1) -- (A3) node [very near start, right=2pt] {$a_{13}$};
  \draw (A2) -- (A4) node [very near end, right=2pt] {$a_{24}$};
  
  \draw [draw=forestgreen, line width=3pt](A1) -- (A2);
  \draw [draw=forestgreen, line width=3pt](A2) -- (A3);
  \draw [draw=forestgreen, line width=3pt](A3) -- (A4);
  \draw [draw=forestgreen, line width=3pt](A4) -- (A1);

  \node [nodeS] (B1) at (8,-1) {\footnotesize 1};
\node [nodeS] (B2) at (11,-1) {\footnotesize 2};
\node [nodeS] (B3) at (11,-4) {\footnotesize 3};
\node [nodeS] (B4) at (8,-4) {\footnotesize 4};
  
  \draw [draw=black, very thick](B1) -- (B2);
  \draw [draw=black, very thick](B1) -- (B3);
  \draw [draw=black, very thick](B2) -- (B3);
  \draw [draw=black, very thick](B2) -- (B4);
  \draw [draw=black, very thick](B3) -- (B4);
  \draw [draw=black, very thick](B4) -- (B1);

  \draw (B1) -- (B2) node [midway, above=2pt] {$a_{12}$};
  \draw (B2) -- (B3) node [midway, right=2pt] {$a_{23}$};
  \draw (B3) -- (B4) node [midway, below=2pt] {$a_{34}$};
  \draw (B4) -- (B1) node [midway, left=2pt] {$a_{14}$};
  \draw (B1) -- (B3) node [very near start, right=2pt] {$a_{13}$};
  \draw (B2) -- (B4) node [very near end, right=2pt] {$a_{24}$};
  
  \draw [draw=blue, line width=3pt](B1) -- (B3);
  \draw [draw=blue, line width=3pt](B2) -- (B3);
  \draw [draw=blue, line width=3pt](B2) -- (B4);
  \draw [draw=blue, line width=3pt](B4) -- (B1);

  \node [nodeS] (C1) at (14,-1) {\footnotesize 1};
\node [nodeS] (C2) at (17,-1) {\footnotesize 2};
\node [nodeS] (C3) at (17,-4) {\footnotesize 3};
\node [nodeS] (C4) at (14,-4) {\footnotesize 4};
  
  \draw [draw=black, very thick](C1) -- (C2);
  \draw [draw=black, very thick](C1) -- (C3);
  \draw [draw=black, very thick](C2) -- (C3);
  \draw [draw=black, very thick](C2) -- (C4);
  \draw [draw=black, very thick](C3) -- (C4);
  \draw [draw=black, very thick](C4) -- (C1);

  \draw (C1) -- (C2) node [midway, above=2pt] {$a_{12}$};
  \draw (C2) -- (C3) node [midway, right=2pt] {$a_{23}$};
  \draw (C3) -- (C4) node [midway, below=2pt] {$a_{34}$};
  \draw (C4) -- (C1) node [midway, left=2pt] {$a_{14}$};
  \draw (C1) -- (C3) node [very near start, right=2pt] {$a_{13}$};
  \draw (C2) -- (C4) node [very near end, right=2pt] {$a_{24}$};
  
  \draw [draw=red, line width=3pt](C1) -- (C3);
  \draw [draw=red, line width=3pt](C1) -- (C2);
  \draw [draw=red, line width=3pt](C2) -- (C4);
  \draw [draw=red, line width=3pt](C3) -- (C4);


\node (T1) at (3.5,-0) { $\boldsymbol{(1,2,3,4,1)}$};
\node (T2) at (9.5,0) {$\boldsymbol{(1,4,2,3,1)}$};
\node (T3) at (15.5,0) {$\boldsymbol{(1,3,4,2,1)}$};

\end{tikzpicture}

    \caption{The three possible $4$-cycles for a $4\times4$ PCM. The elements (above the main diagonal) of the PCM correspond to the edges of the graphs. The three possible $4$-cycles are highlighted by different colors. One can check that Properties~\ref{prop:1}, $\ldots$,~\ref{prop:9} hold.}
    \label{SpecCasesCycles}
\end{figure}

\setcounter{appendbfigure}{2}

\renewcommand{\theproposition}{B\arabic{appendbproposition}}
\begin{proposition}
    A $4\times4$ pairwise comparison matrix $\mathbf{A}$ falls into one of the following categories regarding the number of its consistent $k$-cycles ($k\geq3$):
    \begin{enumerate}
        \item Triple perturbed PCM, i.e., there are no consistent $k$-cycles with $k\geq3$;

        \renewcommand{\labelenumi}{\arabic{enumi}. (a)}
        \item Double perturbed PCM, when there are no consistent $4$-cycles, and exactly one consistent $3$-cycle (triad);

        \renewcommand{\labelenumi}{\arabic{enumi}. (b)}
        \setcounter{enumi}{1}
        \item Double perturbed PCM, when there is exactly one consistent $4$-cycle, and there is no consistent $3$-cycle;

        \renewcommand{\labelenumi}{\arabic{enumi}. (c)}
        \setcounter{enumi}{1}
        \item Double perturbed PCM, when there are exactly two consistent $4$-cycles, and there is no consistent $3$-cycle;

        \renewcommand{\labelenumi}{\arabic{enumi}.}
        \item Simple perturbed PCM, when there is exactly one consistent $4$-cycle, and there are exactly two consistent $3$-cycles;
        \item Consistent PCM, i.e., all $k$-cycles with $k\geq3$ are consistent.
    \end{enumerate}
\end{proposition}
\begin{proof}

    \begin{enumerate}
        \item Triple perturbed case
        
        If there are no consistent $k$-cycles in the $4\times 4$ PCM ($k\geq3$), then it can be made consistent by modifying three elements, as the rank of the matrix determining the triads is three, thus, the PCM is triple perturbed. All other cases must be double perturbed, simple perturbed or consistent.

        \item Double perturbed case

        With the modification of exactly two elements a double perturbed matrix can be made consistent. Based on Properties~\ref{prop:1} and \ref{prop:3}, this modification alters either three or four triads. Three, if the two modified elements appear in a common triad, and four otherwise.

        If all four triads are modified, then each of them is altered only once, thus all of them must be inconsistent at the beginning (and some $4$-cycles must be consistent to get a double perturbed PCM instead of a triple perturbed one).
        
        If three triads are modified, then the fourth one must be consistent at the start. In this case one triad is altered with both elements, and the other two are modified only once, thus these latter two triads must be inconsistent at the beginning. It could be the case that the triad that is modified twice is consistent at the beginning, the modifications reverse each other, and there are two consistent triads. However, based on Property~\ref{prop:8} if there are two consistent triads, then one $4$-cycle must be consistent, and the modification of exactly one element makes the other two $4$-cycles (and thus the PCM) consistent, which means that it is a simple perturbed PCM.

        Thus, for a double perturbed PCM either there is exactly one consistent triad or no consistent triad and one or two consistent $4$-cycles.

        \begin{enumerate}[label=(\alph*)]
            \item If there is exactly one consistent triad, then the other three can be modified and made consistent by altering two elements, which appear in a common triad. Based on Properties~\ref{prop:4} and \ref{prop:9}, however the two elements are chosen this way, it modifies all three $4$-cycles, one of which is modified twice (with both elements), while the other two only altered one time. Thus, the latter two $4$-cycles must be inconsistent before the modification. The last $4$-cycle is altered with both elements, thus, it could be the case that it was consistent at the beginning as well, and the modifications just reverse each other. However, the two modified elements are in a common triad, thus, that would mean that the given triad must have been consistent as well. As there is only one consistent triad, then with these modifications there must be an inconsistent triad that is not modified.

            $\Rightarrow$ If there is exactly one consistent triad in a double perturbed PCM, then all the $4$-cycles must be inconsistent.
            
            \item If there is no consistent triad, then one or two of the $4$-cycles must be consistent to get a double perturbed matrix. If exactly one $4$-cycle is consistent, then based on Properties~\ref{prop:1} and \ref{prop:7}, the two elements that are not included in this $4$-cycle modify all four triads exactly once. Thus, with the modification of those two elements the PCM can be made consistent, this is a possible double perturbed PCM case.
            \item If there is no consistent triad, but two $4$-cycles are consistent that means the PCM is still double perturbed, based on Property~\ref{prop:1}, the four triads can be only modified by at least two elements. According to Properties~\ref{prop:5} and \ref{prop:7}, one can choose any of the two consistent $4$-cycles, and modify those two elements that are not included in the chosen $4$-cycle. In that way one can make all the triads and thus the PCM consistent. Interestingly enough, that means with the first modification one of the two originally consistent $4$-cycles must be made inconsistent again to get a consistent PCM with the next modification. This is also a possible double perturbed PCM case.
        \end{enumerate}

        \item Simple perturbed case

        The PCM can be made consistent by modifying exactly one element. According to Properties~\ref{prop:1} and \ref{prop:2}, the modification of this element will change exactly two triads, and exactly two $4$-cycles, thus this is only possible if there are exactly two consistent triads, and exactly one consistent $4$-cycle in the PCM.

        \item Consistent PCM
        
        By definition, all triads are consistent, and based on Property~\ref{prop:8} those determine all the $4$-cycles as well, which must be consistent.

    \end{enumerate}
\end{proof}

\setcounter{appendbproposition}{2}

\begin{proposition} \label{setofefficientweightvectorsspeccases}
    The set of efficient weight vectors of a $4\times4$ pairwise comparison matrix $\mathbf{A}$ is the union of three tetrahedra (see Proposition~\ref{final} of the main text of the paper) with the following special cases:
    \begin{enumerate}
        \item For a triple perturbed PCM, the three tetrahedra have no common vertex, however,  for each pair of tetrahedra, for the two tetrahedra in the given pair, one of their edges is on the same line, and two of their faces are pairwise on the same plane. (see Example~\ref{mainexample} of the main text of the paper);

        \renewcommand{\labelenumi}{\arabic{enumi}. (a)}
        \item For a double perturbed PCM, when there are no consistent $4$-cycles, and exactly one consistent $3$-cycle (triad), each pair of tetrahedra have a common vertex  (see Example~\ref{doubletriad} as well as \citet[Figure 6]{Szadoczki2024});

        \renewcommand{\labelenumi}{\arabic{enumi}. (b)}
        \setcounter{enumi}{1}
        \item For a double perturbed PCM, when there are exactly one consistent $4$-cycle, and there is no consistent $3$-cycle, one tetrahedron corresponds to a point, and the other two do not have a vertex which corresponds to this single point, nor a common vertex with each other (see Example~\ref{doubleoneconsistentfourcycle});

        \renewcommand{\labelenumi}{\arabic{enumi}. (c)}
        \setcounter{enumi}{1}
        \item For a double perturbed PCM, when there are exactly two consistent $4$-cycles, and there is no consistent $3$-cycle, two of the tetrahedra correspond to a different point each that is not common with any of the vertices of the third tetrahedron (see Example~\ref{doubletwoconsistentfourcycles});

        \renewcommand{\labelenumi}{\arabic{enumi}.}
        \item For a simple perturbed PCM, when there is exactly one consistent $4$-cycle, and there are exactly two consistent $3$-cycles, one of the tetrahedra corresponds to a point that is the same as one of the vertices of the other two tetrahedra. For those two tetrahedra, three of their four vertices pairwise correspond, i.e., they have a common face (see Example~\ref{simpleperturbed} as well as \citet[Figure 8]{Szadoczki2024});
        \item For a consistent PCM, the set of efficient weight vectors corresponds to a single vector, a single point (see Example~\ref{consistent}).
    \end{enumerate}
\end{proposition}

\begin{proof}
The set of efficient weight vectors is connected, and Proposition~\ref{final} holds without Assumption~\ref{assumption:1} for all cases as well. The only difference is that, it can happen that more than three directed Hamiltonian cycles are possible, but that means there are consistent $4$-cycles in the PCM, and all spanning trees obtained from those will determine the same weight vector.

Thus, only the connection of tetrahedra (i.e., the correspondence of vertices) must be proven. The different vertices (different weight vectors calculated from different incomplete pairwise comparison matrices of the PCM with path spanning tree representing graphs) can only correspond if there are consistent $k$-cycles ($k\geq3$) in the PCM.

\begin{enumerate}
    \item Triple perturbed PCM

    There are no consistent $k$-cycles ($k\geq3$), thus, the three tetrahedra have no common vertex, Proposition~\ref{final} holds with Assumption~\ref{assumption:1}.
    
    Based on Property~\ref{prop:5} there are two common elements for each pair of the possible $4$-cycles. The path spanning trees contain three elements of the four for a given $4$-cycle. That means for a given pair of tetrahedra (determined by two $4$-cycles), for each tetrahedron in the pair, two of their four spanning trees have two common elements with the other $4$-cycle in the pair (the calculated weight vectors are on the same line), and there are two triplets of spanning trees for each tetrahedron that have one common element with the other $4$-cycle in the pair (the calculated weight vectors are on the same cutting plane). Thus, for the two tetrahedra one of their edges is on the same line, and two of their faces are pairwise on the same plane.

    \renewcommand{\labelenumi}{\arabic{enumi}. (a)}
    \item Double perturbed PCM, when there are no consistent $4$-cycles, and exactly one consistent $3$-cycle (triad)

     Based on Properties~\ref{prop:4}, \ref{prop:5}, and \ref{prop:6}, for each pair of tetrahedra there is a path spanning tree representing graph for each of them that contains two elements from the consistent triad, and the third element is the same. As the consistent triad contains the same information, the weight vectors (and thus the appropriate vertices) will correspond. That means each pair of tetrahedra have exactly one common vertex.

    \renewcommand{\labelenumi}{\arabic{enumi}. (b)}
        \setcounter{enumi}{1}
        \item Double perturbed PCM, when there are exactly one consistent $4$-cycle, and there is no consistent $3$-cycle

        All the weight vectors calculated from the path spanning trees of the consistent $4$-cycle will correspond. Thus the given tetrahedron corresponds to a single point. There are no further consistent $4$-cycles or triads, thus, no further path spanning tree will provide the same weight vector. 

        \renewcommand{\labelenumi}{\arabic{enumi}. (c)}
        \setcounter{enumi}{1}
        \item Double perturbed PCM, when there are exactly two consistent $4$-cycles, and there is no consistent $3$-cycle

        All the weight vectors calculated from the path spanning trees of each of the consistent $4$-cycles will correspond. However, according to Properties~\ref{prop:1}, \ref{prop:5} and \ref{prop:7}, the vectors that are calculated from different consistent $4$-cycles cannot correspond to each other, because that would mean that every element is consistent with each other in the PCM. As there are no further consistent $k$-cycles ($k\geq3$) in the PCM, no further vertices will correspond to each other.
        
        That means two of the tetrahedra that are coming from the consistent $4$-cycles correspond to a different point each that is not common with any of the vertices of the third tetrahedron.

        \renewcommand{\labelenumi}{\arabic{enumi}.}
        \item Simple perturbed PCM, when there is exactly one consistent $4$-cycle, and there are exactly two consistent $3$-cycles

       The weight vectors calculated from the path spanning trees of the consistent $4$-cycle will correspond, thus one tetrahedron will correspond to a single point. Based on Property~\ref{prop:8}, as there are exactly two consistent triads, those must determine the only consistent $4$-cycle, and there is exactly one element (let us denote it by $a_{sp}$) that is not included in their union, but it is included in the other two possible $4$-cycles. Modifying $a_{sp}$ makes the PCM consistent. However, that means for each of the inconsistent $4$-cycles there is a path representing spanning tree that excludes $a_{sp}$, and this spanning tree is consistent with the consistent $4$-cycle, thus it produces the same weight vector (vertex). This means that the tetrahedron that corresponds to a single point will be also a vertex of both of the other tetrahedra.

       Based on Property~\ref{prop:4}, for each consistent triad there will be exactly one path spanning tree for each $4$-cycles that contains two elements from the given triad, and $a_{sp}$. Those will also provide the same weight vectors, as their inconsistency coming from the same element, but a different one than the one calculated from the consistent $4$-cycle (as there is inconsistency compared to that). That means for the two tetrahedra coming from the inconsistent $4$-cycles three of their four vertices pairwise correspond, and one pair also corresponds with the third (consistent) tetrahedron.

       Based on Properties~\ref{prop:5}, \ref{prop:6}, and \ref{prop:7}, for the two inconsistent $4$-cycles one of their two common elements is $a_{sp}$, and the other must be the part of both consistent triads. This means that there is a path spanning tree for each that includes $a_{sp}$, and two other elements that are not in the other inconsistent $4$-cycle, and they are also not in the same consistent triad. These spanning trees will result in different weight vectors from each other, and different ones from the so far discussed vertices as well.

        Thus, the remaining two vertices differ from the previous ones and differ from each other as well in the appropriate tetrahedra.

        \item Consistent PCM
        
        The entries of a consistent PCM $\mathbf{A}$ can be written as the ratios of the elements of exactly one (normalized) weight vector $\mathbf{w}$, thus the only efficient weight vector in this case is $\mathbf{w}$.
\end{enumerate}
\end{proof}

We show an example for each special case modifying the triple perturbed PCM used in Example~\ref{mainexample}:

\[ \mathbf{A''} =
\begin{pmatrix}
1 & 1 & 5 & 7 \\
1 & 1 & 2 & 8 \\
1/5 & 1/2 & 1 & 1/3 \\
1/7 & 1/8 & 3 & 1
\end{pmatrix}.
\]

\renewcommand{\theexample}{B\arabic{appendbexample}}
\begin{example} \label{doubletriad}
Figure~\ref{DoublePerturbed_triad} shows the set of efficient weight vectors for the following double perturbed PCM, which has exactly one consistent triad, and no consistent $4$-cycle:

\[ \mathbf{A''} =
\begin{pmatrix}
1 & \color{orange}5/2 & 5 & 7 \\
\color{orange}2/5 & 1 & 2 & 8 \\
1/5 & 1/2 & 1 & 1/3 \\
1/7 & 1/8 & 3 & 1
\end{pmatrix}.
\]

Here the original example is only modified by one element (and its reciprocal), which is highlighted by orange color in the matrix. All of the $4$-cycles of the matrix are inconsistent, however, the $(1,2,3)$ triad is consistent, and the vertices calculated from two elements of this triad, and the same third element correspond to each other. That means each pair of tetrahedra has one common vertex in Figure~\ref{DoublePerturbed_triad}.

\renewcommand{\thefigure}{B\arabic{appendbfigure}}
\begin{figure}[H]
    \centering

\begin{tikzpicture}[line join = round, line cap = round, scale=0.87]

\pgfmathsetmacro{\mainfact}{8};

\coordinate [label=above left:\text{$V_1=(w_1=1,w_2=0,w_3=0,w_4=0)$}] (A) at (1*\mainfact,1*\mainfact,0);
\coordinate [label=below:\text{$V_2=(w_1=0,w_2=1,w_3=0,w_4=0)$}] (B) at (1*\mainfact,0,1*\mainfact);
\coordinate [label=left:\text{$V_3=(w_1=0,w_2=0,$}] (C) at (0,1*\mainfact,1*\mainfact);
\coordinate [label=left:\text{$w_3=1,w_4=0)$}] (CC) at (0,0.93*\mainfact,1*\mainfact);
\coordinate [label=left:\text{$V_4=(w_1=0,w_2=0,w_3=0,w_4=1)$}] (D) at (0,0,0);

\draw[->] (0,0) -- (10,0,0) node[right] {$x$};
\draw[->] (0,0) -- (0,10,0) node[above] {$y$};
\draw[->] (0,0) -- (0,0,10) node[below left] {$z$};

\foreach \i in {A,B,C,D}
    \draw[dashed,very thick] (0,0)--(\i);
\draw[-,very thick, opacity=.5] (A)--(D)--(B)--cycle;
\draw[-,very thick, opacity=.5] (A) --(D)--(C)--cycle;
\draw[-,very thick, opacity=.5] (B)--(D)--(C)--cycle;
\draw[-,very thick] (A)--(B);
\draw[-,very thick] (A)--(C);
\draw[-,very thick] (A)--(D);

\coordinate (M) at (0.954545455*\mainfact,0.693181818*\mainfact,0.284090909*\mainfact);
\coordinate (N) at (0.848484848*\mainfact,0.727272727*\mainfact,0.363636364*\mainfact);
\coordinate (O) at (0.636363636*\mainfact,0.545454545*\mainfact,0.272727273*\mainfact);
\coordinate (P) at (0.878787879*\mainfact,0.181818182*\mainfact,0.757575758*\mainfact);

\draw[-,very thick, fill = red!30, opacity =.5] (M)--(N)--(O)--cycle;
\draw[-,very thick, fill = red!30, opacity =.5] (M) --(O)--(P)--cycle;
\draw[-,very thick, fill = red!30, opacity =.5] (N)--(O)--(P)--cycle;

\coordinate (I) at (0.848484848*\mainfact,0.727272727*\mainfact,0.363636364*\mainfact);
\coordinate (J) at (0.803278689*\mainfact,0.68852459*\mainfact,0.344262295*\mainfact);
\coordinate (K) at (0.862068966*\mainfact,0.482758621*\mainfact,0.540229885*\mainfact);
\coordinate (L) at (0.75*\mainfact,0.55*\mainfact,0.6*\mainfact);

\draw[-,very thick, fill = blue!30, opacity =.5] (I)--(J)--(K)--cycle;
\draw[-,very thick, fill = blue!30, opacity =.5] (I) --(K)--(L)--cycle;
\draw[-,very thick, fill = blue!30, opacity =.5] (J)--(K)--(L)--cycle;

\coordinate (E) at (0.636363636*\mainfact,0.545454545*\mainfact,0.272727273*\mainfact);
\coordinate (F) at (0.803278689*\mainfact,0.68852459*\mainfact,0.344262295*\mainfact);
\coordinate (G) at (0.880239521*\mainfact,0.658682635*\mainfact,0.281437126*\mainfact);
\coordinate (H) at (0.851851852*\mainfact,0.814814815*\mainfact,0.111111111*\mainfact);

\draw[-,very thick, fill=green!30, opacity=.5] (E)--(H)--(F)--cycle;
\draw[-,very thick, fill=green!30, opacity=.5] (E) --(H)--(G)--cycle;
\draw[-,very thick, fill=green!30, opacity=.5] (F)--(H)--(G)--cycle;

\tikzstyle{node1} = [circle,inner sep=2pt,draw=forestgreen,fill=forestgreen,thick];

\node [node1] (A1) at (8.5,7.5) {\footnotesize \color{white} 1};
\node [node1] (A2) at (10,7.5) {\footnotesize \color{white}  2};
\node [node1] (A3) at (10,6) {\footnotesize \color{white} 3};
\node [node1] (A4) at (8.5,6) {\footnotesize \color{white} 4};

  \draw [-latex,draw=forestgreen, very thick](A1) -- (A2);
  \draw [-latex,draw=forestgreen, very thick](A2) -- (A3);
  \draw [-latex,draw=forestgreen, very thick](A3) -- (A4);
  \draw [-latex,draw=forestgreen, very thick](A4) -- (A1);

  \tikzstyle{node2} = [circle,inner sep=2pt,draw=blue,fill= blue,thick];

\node [node2] (B1) at (8.5,5) {\footnotesize \color{white} 1};
\node [node2] (B2) at (8.5,3.5) {\footnotesize \color{white} 3};
\node [node2] (B3) at (10,3.5) {\footnotesize \color{white} 2};
\node [node2] (B4) at (10,5) {\footnotesize \color{white} 4};

  \draw [-latex,draw=blue, very thick](B2) -- (B1);
  \draw [-latex,draw=blue, very thick](B3) -- (B2);
  \draw [-latex,draw=blue, very thick](B4) -- (B3);
  \draw [-latex,draw=blue, very thick](B1) -- (B4);

\tikzstyle{node3} = [circle,inner sep=2pt,draw=red,fill= red,thick];

\node [node3] (C1) at (8.5,2.5) {\footnotesize \color{white} 1};
\node [node3] (C2) at (8.5,1) {\footnotesize \color{white} 2};
\node [node3] (C3) at (10,1) {\footnotesize \color{white} 4};
\node [node3] (C4) at (10,2.5) {\footnotesize \color{white} 3};

  \draw [-latex,draw=red, very thick](C2) -- (C1);
  \draw [-latex,draw=red, very thick](C3) -- (C2);
  \draw [-latex,draw=red, very thick](C4) -- (C3);
  \draw [-latex,draw=red, very thick](C1) -- (C4);

\end{tikzpicture}

    \caption{The set of efficient weight vectors (the union of $\mathcal{T}_1, \mathcal{T}_2$ and $\mathcal{T}_3$) for the double perturbed example with exactly one (the $(1,2,3)$) consistent triad. An interactive version is available at https://www.geogebra.org/m/jw8xu5ud.}
    \label{DoublePerturbed_triad}
\end{figure}

\setcounter{appendbfigure}{3}

\end{example}

\setcounter{appendbexample}{2}

\begin{example} \label{doubleoneconsistentfourcycle}
Figure~\ref{DoublePerturbed} shows the set of efficient weight vectors for the following double perturbed PCM, which has no consistent triad, and exactly one consistent $4$-cycle:

\[ \mathbf{A''} =
\begin{pmatrix}
1 & 1 & 5 & 7 \\
1 & 1 & 2 & \color{blue} 14/5 \color{black} \\
1/5 & 1/2 & 1 & 1/3 \\
1/7 & \color{blue}5/14 \color{black} & 3 & 1
\end{pmatrix}.
\]

Here the original example is only modified by one element (and its reciprocal), which is highlighted by blue color in the matrix. All of the triads of the matrix are inconsistent, however, the $(1,4,2,3,1)$ $4$-cycle is consistent, and the appropriate tetrahedron corresponds to one point in Figure~\ref{DoublePerturbed}.

\renewcommand{\thefigure}{B\arabic{appendbfigure}}
\begin{figure}[H]
    \centering

\begin{tikzpicture}[line join = round, line cap = round, scale=0.87]

\pgfmathsetmacro{\mainfact}{8};

\coordinate [label=above left:\text{$V_1=(w_1=1,w_2=0,w_3=0,w_4=0)$}] (A) at (1*\mainfact,1*\mainfact,0);
\coordinate [label=below:\text{$V_2=(w_1=0,w_2=1,w_3=0,w_4=0)$}] (B) at (1*\mainfact,0,1*\mainfact);
\coordinate [label=left:\text{$V_3=(w_1=0,w_2=0,$}] (C) at (0,1*\mainfact,1*\mainfact);
\coordinate [label=left:\text{$w_3=1,w_4=0)$}] (CC) at (0,0.93*\mainfact,1*\mainfact);
\coordinate [label=left:\text{$V_4=(w_1=0,w_2=0,w_3=0,w_4=1)$}] (D) at (0,0,0);

\draw[->] (0,0) -- (10,0,0) node[right] {$x$};
\draw[->] (0,0) -- (0,10,0) node[above] {$y$};
\draw[->] (0,0) -- (0,0,10) node[below left] {$z$};

\foreach \i in {A,B,C,D}
    \draw[dashed,very thick] (0,0)--(\i);
\draw[-,very thick, opacity=.5] (A)--(D)--(B)--cycle;
\draw[-,very thick, opacity=.5] (A) --(D)--(C)--cycle;
\draw[-,very thick, opacity=.5] (B)--(D)--(C)--cycle;
\draw[-,very thick] (A)--(B);
\draw[-,very thick] (A)--(C);
\draw[-,very thick] (A)--(D);

\coordinate (E) at (0.5*\mainfact,0.375*\mainfact,0.375*\mainfact);
\coordinate (F) at (0.756756757*\mainfact,0.567567568*\mainfact,0.567567568*\mainfact);
\coordinate (G) at (0.913043478*\mainfact,0.47826087*\mainfact,0.47826087*\mainfact);
\coordinate (H) at (0.851851852*\mainfact,0.814814815*\mainfact,0.111111111*\mainfact);

\draw[-,very thick, fill=green!30, opacity=.5] (E)--(H)--(F)--cycle;
\draw[-,very thick, fill=green!30, opacity=.5] (E) --(H)--(G)--cycle;
\draw[-,very thick, fill=green!30, opacity=.5] (F)--(H)--(G)--cycle;

\coordinate (I) at (0.803278689*\mainfact,0.68852459*\mainfact,0.344262295*\mainfact);
\coordinate (J) at (0.803278689*\mainfact,0.68852459*\mainfact,0.344262295*\mainfact);
\coordinate (K) at (0.803278689*\mainfact,0.68852459*\mainfact,0.344262295*\mainfact);
\coordinate (L) at (0.803278689*\mainfact,0.68852459*\mainfact,0.344262295*\mainfact);

\draw[-,very thick, fill = blue!30, opacity =.5] (I)--(J)--(K)--cycle;
\draw[-,very thick, fill = blue!30, opacity =.5] (I) --(K)--(L)--cycle;
\draw[-,very thick, fill = blue!30, opacity =.5] (J)--(K)--(L)--cycle;

\coordinate (M) at (0.807692308*\mainfact,0.451923077*\mainfact,0.451923077*\mainfact);
\coordinate (N) at (0.782122905*\mainfact,0.469273743*\mainfact,0.469273743*\mainfact);
\coordinate (O) at (0.714285714*\mainfact,0.428571429*\mainfact,0.428571429*\mainfact);
\coordinate (P) at (0.770114943*\mainfact,0.344827586*\mainfact,0.540229885*\mainfact);

\draw[-,very thick, fill = red!30, opacity =.5] (M)--(N)--(O)--cycle;
\draw[-,very thick, fill = red!30, opacity =.5] (M) --(O)--(P)--cycle;
\draw[-,very thick, fill = red!30, opacity =.5] (N)--(O)--(P)--cycle;

\tikzstyle{node1} = [circle,inner sep=2pt,draw=forestgreen,fill=forestgreen,thick];

\node [node1] (A1) at (8.5,7.5) {\footnotesize \color{white} 1};
\node [node1] (A2) at (10,7.5) {\footnotesize \color{white}  2};
\node [node1] (A3) at (10,6) {\footnotesize \color{white} 3};
\node [node1] (A4) at (8.5,6) {\footnotesize \color{white} 4};

  \draw [-latex,draw=forestgreen, very thick](A1) -- (A2);
  \draw [-latex,draw=forestgreen, very thick](A2) -- (A3);
  \draw [-latex,draw=forestgreen, very thick](A3) -- (A4);
  \draw [-latex,draw=forestgreen, very thick](A4) -- (A1);

  \tikzstyle{node2} = [circle,inner sep=2pt,draw=blue,fill= blue,thick];

\node [node2] (B1) at (8.5,5) {\footnotesize \color{white} 1};
\node [node2] (B2) at (8.5,3.5) {\footnotesize \color{white} 3};
\node [node2] (B3) at (10,3.5) {\footnotesize \color{white} 2};
\node [node2] (B4) at (10,5) {\footnotesize \color{white} 4};

  \draw [latex-latex,draw=blue, very thick](B2) -- (B1);
  \draw [latex-latex,draw=blue, very thick](B3) -- (B2);
  \draw [latex-latex,draw=blue, very thick](B4) -- (B3);
  \draw [latex-latex,draw=blue, very thick](B1) -- (B4);

\tikzstyle{node3} = [circle,inner sep=2pt,draw=red,fill= red,thick];

\node [node3] (C1) at (8.5,2.5) {\footnotesize \color{white} 1};
\node [node3] (C2) at (8.5,1) {\footnotesize \color{white} 2};
\node [node3] (C3) at (10,1) {\footnotesize \color{white} 4};
\node [node3] (C4) at (10,2.5) {\footnotesize \color{white} 3};

  \draw [-latex,draw=red, very thick](C2) -- (C1);
  \draw [-latex,draw=red, very thick](C3) -- (C2);
  \draw [-latex,draw=red, very thick](C4) -- (C3);
  \draw [-latex,draw=red, very thick](C1) -- (C4);

  \tikzstyle{node4} = [circle,inner sep=1pt,draw=blue,fill=blue,opacity=0.5];

\node [node4] (C1423) at (0.803278689*\mainfact,0.68852459*\mainfact,0.344262295*\mainfact) {};

\end{tikzpicture}

    \caption{The set of efficient weight vectors (the union of $\mathcal{T}_1, \mathcal{T}_2$ and $\mathcal{T}_3$) for the double perturbed example with no consistent triad, and exactly one (the $(1,4,2,3,1)$) consistent $4$-cycle. An interactive version is available at https://www.geogebra.org/m/sxdnxxru.}
    \label{DoublePerturbed}
\end{figure}

\setcounter{appendbfigure}{4}

\end{example}

\setcounter{appendbexample}{3}
\begin{example} \label{doubletwoconsistentfourcycles}
Figure~\ref{DoublePerturbed2consistent4cycles} shows the set of efficient weight vectors for the following double perturbed PCM, which has no consistent triad, and exactly two consistent $4$-cycles:

\[ \mathbf{A''} =
\begin{pmatrix}
1 & 1 & 5 & 7 \\
1 & 1 & 2 & \color{blue} 14/5 \color{black} \\
1/5 & 1/2 & 1 & \color{red}14/25 \color{black} \\
1/7 & \color{blue}5/14 \color{black} & \color{red}25/14\color{black} & 1
\end{pmatrix}.
\]

Here the original example is only modified by two elements (and their reciprocals), which are highlighted by blue and red colors in the matrix. All of the triads of the matrix are inconsistent, however, the $(1,4,2,3,1)$ and the $(1,3,4,2,1)$ $4$-cycles are consistent, and the appropriate tetrahedra correspond to two different points in Figure~\ref{DoublePerturbed2consistent4cycles}.

\renewcommand{\thefigure}{B\arabic{appendbfigure}}
\begin{figure}[H]
    \centering

\begin{tikzpicture}[line join = round, line cap = round, scale=0.87]

\pgfmathsetmacro{\mainfact}{8};

\coordinate [label=above left:\text{$V_1=(w_1=1,w_2=0,w_3=0,w_4=0)$}] (A) at (1*\mainfact,1*\mainfact,0);
\coordinate [label=below:\text{$V_2=(w_1=0,w_2=1,w_3=0,w_4=0)$}] (B) at (1*\mainfact,0,1*\mainfact);
\coordinate [label=left:\text{$V_3=(w_1=0,w_2=0,$}] (C) at (0,1*\mainfact,1*\mainfact);
\coordinate [label=left:\text{$w_3=1,w_4=0)$}] (CC) at (0,0.93*\mainfact,1*\mainfact);
\coordinate [label=left:\text{$V_4=(w_1=0,w_2=0,w_3=0,w_4=1)$}] (D) at (0,0,0);

\draw[->] (0,0) -- (10,0,0) node[right] {$x$};
\draw[->] (0,0) -- (0,10,0) node[above] {$y$};
\draw[->] (0,0) -- (0,0,10) node[below left] {$z$};

\foreach \i in {A,B,C,D}
    \draw[dashed,very thick] (0,0)--(\i);
\draw[-,very thick, opacity=.5] (A)--(D)--(B)--cycle;
\draw[-,very thick, opacity=.5] (A) --(D)--(C)--cycle;
\draw[-,very thick, opacity=.5] (B)--(D)--(C)--cycle;
\draw[-,very thick] (A)--(B);
\draw[-,very thick] (A)--(C);
\draw[-,very thick] (A)--(D);

\coordinate (E) at (0.589473684*\mainfact,0.442105263*\mainfact,0.442105263*\mainfact);
\coordinate (F) at (0.756756757*\mainfact,0.567567568*\mainfact,0.567567568*\mainfact);
\coordinate (G) at (0.899742931*\mainfact,0.485861183*\mainfact,0.485861183*\mainfact);
\coordinate (H) at (0.838842975*\mainfact,0.780991736*\mainfact,0.173553719*\mainfact);

\draw[-,very thick, fill=green!30, opacity=.5] (E)--(H)--(F)--cycle;
\draw[-,very thick, fill=green!30, opacity=.5] (E) --(H)--(G)--cycle;
\draw[-,very thick, fill=green!30, opacity=.5] (F)--(H)--(G)--cycle;

\coordinate (I) at (0.803278689*\mainfact,0.68852459*\mainfact,0.344262295*\mainfact);
\coordinate (J) at (0.803278689*\mainfact,0.68852459*\mainfact,0.344262295*\mainfact);
\coordinate (K) at (0.803278689*\mainfact,0.68852459*\mainfact,0.344262295*\mainfact);
\coordinate (L) at (0.803278689*\mainfact,0.68852459*\mainfact,0.344262295*\mainfact);

\draw[-,very thick, fill = blue!30, opacity =.5] (I)--(J)--(K)--cycle;
\draw[-,very thick, fill = blue!30, opacity =.5] (I) --(K)--(L)--cycle;
\draw[-,very thick, fill = blue!30, opacity =.5] (J)--(K)--(L)--cycle;

\coordinate (M) at (0.782122905*\mainfact,0.469273743*\mainfact,0.469273743*\mainfact);
\coordinate (N) at (0.782122905*\mainfact,0.469273743*\mainfact,0.469273743*\mainfact);
\coordinate (O) at (0.782122905*\mainfact,0.469273743*\mainfact,0.469273743*\mainfact);
\coordinate (P) at (0.782122905*\mainfact,0.469273743*\mainfact,0.469273743*\mainfact);

\draw[-,very thick, fill = red!30, opacity =.5] (M)--(N)--(O)--cycle;
\draw[-,very thick, fill = red!30, opacity =.5] (M) --(O)--(P)--cycle;
\draw[-,very thick, fill = red!30, opacity =.5] (N)--(O)--(P)--cycle;

\tikzstyle{node1} = [circle,inner sep=2pt,draw=forestgreen,fill=forestgreen,thick];

\node [node1] (A1) at (8.5,7.5) {\footnotesize \color{white} 1};
\node [node1] (A2) at (10,7.5) {\footnotesize \color{white}  2};
\node [node1] (A3) at (10,6) {\footnotesize \color{white} 3};
\node [node1] (A4) at (8.5,6) {\footnotesize \color{white} 4};

  \draw [-latex,draw=forestgreen, very thick](A1) -- (A2);
  \draw [-latex,draw=forestgreen, very thick](A2) -- (A3);
  \draw [-latex,draw=forestgreen, very thick](A3) -- (A4);
  \draw [-latex,draw=forestgreen, very thick](A4) -- (A1);

  \tikzstyle{node2} = [circle,inner sep=2pt,draw=blue,fill= blue,thick];

\node [node2] (B1) at (8.5,5) {\footnotesize \color{white} 1};
\node [node2] (B2) at (8.5,3.5) {\footnotesize \color{white} 3};
\node [node2] (B3) at (10,3.5) {\footnotesize \color{white} 2};
\node [node2] (B4) at (10,5) {\footnotesize \color{white} 4};

  \draw [latex-latex,draw=blue, very thick](B2) -- (B1);
  \draw [latex-latex,draw=blue, very thick](B3) -- (B2);
  \draw [latex-latex,draw=blue, very thick](B4) -- (B3);
  \draw [latex-latex,draw=blue, very thick](B1) -- (B4);

\tikzstyle{node3} = [circle,inner sep=2pt,draw=red,fill= red,thick];

\node [node3] (C1) at (8.5,2.5) {\footnotesize \color{white} 1};
\node [node3] (C2) at (8.5,1) {\footnotesize \color{white} 2};
\node [node3] (C3) at (10,1) {\footnotesize \color{white} 4};
\node [node3] (C4) at (10,2.5) {\footnotesize \color{white} 3};

  \draw [latex-latex,draw=red, very thick](C2) -- (C1);
  \draw [latex-latex,draw=red, very thick](C3) -- (C2);
  \draw [latex-latex,draw=red, very thick](C4) -- (C3);
  \draw [latex-latex,draw=red, very thick](C1) -- (C4);

  \tikzstyle{node4} = [circle,inner sep=1pt,draw=blue,fill=blue,opacity=0.5];

  \tikzstyle{node5} = [circle,inner sep=1pt,draw=red,fill=red,opacity=0.5];

\node [node4] (C1423) at (0.803278689*\mainfact,0.68852459*\mainfact,0.344262295*\mainfact) {};

\node [node5] (C1342) at (0.782122905*\mainfact,0.469273743*\mainfact,0.469273743*\mainfact) {};

\end{tikzpicture}

    \caption{The set of efficient weight vectors (the union of $\mathcal{T}_1, \mathcal{T}_2$ and $\mathcal{T}_3$) for the double perturbed example with no consistent triad, and exactly two (the $(1,4,2,3,1)$ and $(1,3,4,2,1)$) consistent $4$-cycles.}
    \label{DoublePerturbed2consistent4cycles}
\end{figure}

\setcounter{appendbfigure}{5}
\end{example}

\setcounter{appendbexample}{4}
\begin{example} \label{simpleperturbed}
Figure~\ref{SimplePerturbed} shows the set of efficient weight vectors for the following simple perturbed PCM:

\[ \mathbf{A''} =
\begin{pmatrix}
1 & \color{orange} 5/2 & 5 & 7 \\
\color{orange}2/5 & 1 & 2 & \color{blue} 14/5 \color{black} \\
1/5 & 1/2 & 1 & 1/3 \\
1/7 & \color{blue}5/14 \color{black} & 3 & 1
\end{pmatrix}.
\]

Here the original example is only modified by two elements (and their reciprocals), which are highlighted by orange and blue colors in the matrix. The $(1,2,3)$ and $(1,2,4)$ triads, and the $(1,4,2,3,1)$ $4$-cycle are consistent. The appropriate tetrahedron corresponds to a point that is a vertex of the other two tetrahedra. For those two tetrahedra, three of their four vertices pairwise correspond in  Figure~\ref{SimplePerturbed}.

\renewcommand{\thefigure}{B\arabic{appendbfigure}}
\begin{figure}[H]
    \centering

\begin{tikzpicture}[line join = round, line cap = round, scale=0.87]

\pgfmathsetmacro{\mainfact}{8};

\coordinate [label=above left:\text{$V_1=(w_1=1,w_2=0,w_3=0,w_4=0)$}] (A) at (1*\mainfact,1*\mainfact,0);
\coordinate [label=below:\text{$V_2=(w_1=0,w_2=1,w_3=0,w_4=0)$}] (B) at (1*\mainfact,0,1*\mainfact);
\coordinate [label=left:\text{$V_3=(w_1=0,w_2=0,$}] (C) at (0,1*\mainfact,1*\mainfact);
\coordinate [label=left:\text{$w_3=1,w_4=0)$}] (CC) at (0,0.93*\mainfact,1*\mainfact);
\coordinate [label=left:\text{$V_4=(w_1=0,w_2=0,w_3=0,w_4=1)$}] (D) at (0,0,0);

\draw[->] (0,0) -- (10,0,0) node[right] {$x$};
\draw[->] (0,0) -- (0,10,0) node[above] {$y$};
\draw[->] (0,0) -- (0,0,10) node[below left] {$z$};

\foreach \i in {A,B,C,D}
    \draw[dashed,very thick] (0,0)--(\i);
\draw[-,very thick, opacity=.5] (A)--(D)--(B)--cycle;
\draw[-,very thick, opacity=.5] (A) --(D)--(C)--cycle;
\draw[-,very thick, opacity=.5] (B)--(D)--(C)--cycle;
\draw[-,very thick] (A)--(B);
\draw[-,very thick] (A)--(C);
\draw[-,very thick] (A)--(D);

\coordinate (I) at (0.803278689*\mainfact,0.68852459*\mainfact,0.344262295*\mainfact);
\coordinate (J) at (0.803278689*\mainfact,0.68852459*\mainfact,0.344262295*\mainfact);
\coordinate (K) at (0.803278689*\mainfact,0.68852459*\mainfact,0.344262295*\mainfact);
\coordinate (L) at (0.803278689*\mainfact,0.68852459*\mainfact,0.344262295*\mainfact);

\draw[-,very thick, fill = blue!30, opacity =.5] (I)--(J)--(K)--cycle;
\draw[-,very thick, fill = blue!30, opacity =.5] (I) --(K)--(L)--cycle;
\draw[-,very thick, fill = blue!30, opacity =.5] (J)--(K)--(L)--cycle;

\coordinate (M) at (0.880239521*\mainfact,0.658682635*\mainfact,0.281437126*\mainfact);
\coordinate (N) at (0.803278689*\mainfact,0.68852459*\mainfact,0.344262295*\mainfact);
\coordinate (O) at (0.636363636*\mainfact,0.545454545*\mainfact,0.272727273*\mainfact);
\coordinate (P) at (0.770114943*\mainfact,0.344827586*\mainfact,0.540229885*\mainfact);

\draw[-,very thick, fill = red!30, opacity =.5] (M)--(N)--(O)--cycle;
\draw[-,very thick, fill = red!30, opacity =.5] (M) --(O)--(P)--cycle;
\draw[-,very thick, fill = red!30, opacity =.5] (N)--(O)--(P)--cycle;

\coordinate (E) at (0.636363636*\mainfact,0.545454545*\mainfact,0.272727273*\mainfact);
\coordinate (F) at (0.803278689*\mainfact,0.68852459*\mainfact,0.344262295*\mainfact);
\coordinate (G) at (0.880239521*\mainfact,0.658682635*\mainfact,0.281437126*\mainfact);
\coordinate (H) at (0.851851852*\mainfact,0.814814815*\mainfact,0.111111111*\mainfact);

\draw[-,very thick, fill=green!30, opacity=.5] (E)--(H)--(F)--cycle;
\draw[-,very thick, fill=green!30, opacity=.5] (E) --(H)--(G)--cycle;
\draw[-,very thick, fill=green!30, opacity=.5] (F)--(H)--(G)--cycle;

\tikzstyle{node1} = [circle,inner sep=2pt,draw=forestgreen,fill=forestgreen,thick];

\node [node1] (A1) at (8.5,7.5) {\footnotesize \color{white} 1};
\node [node1] (A2) at (10,7.5) {\footnotesize \color{white}  2};
\node [node1] (A3) at (10,6) {\footnotesize \color{white} 3};
\node [node1] (A4) at (8.5,6) {\footnotesize \color{white} 4};

  \draw [-latex,draw=forestgreen, very thick](A1) -- (A2);
  \draw [-latex,draw=forestgreen, very thick](A2) -- (A3);
  \draw [-latex,draw=forestgreen, very thick](A3) -- (A4);
  \draw [-latex,draw=forestgreen, very thick](A4) -- (A1);

  \tikzstyle{node2} = [circle,inner sep=2pt,draw=blue,fill= blue,thick];

\node [node2] (B1) at (8.5,5) {\footnotesize \color{white} 1};
\node [node2] (B2) at (8.5,3.5) {\footnotesize \color{white} 3};
\node [node2] (B3) at (10,3.5) {\footnotesize \color{white} 2};
\node [node2] (B4) at (10,5) {\footnotesize \color{white} 4};

  \draw [latex-latex,draw=blue, very thick](B2) -- (B1);
  \draw [latex-latex,draw=blue, very thick](B3) -- (B2);
  \draw [latex-latex,draw=blue, very thick](B4) -- (B3);
  \draw [latex-latex,draw=blue, very thick](B1) -- (B4);

\tikzstyle{node3} = [circle,inner sep=2pt,draw=red,fill= red,thick];

\node [node3] (C1) at (8.5,2.5) {\footnotesize \color{white} 1};
\node [node3] (C2) at (8.5,1) {\footnotesize \color{white} 2};
\node [node3] (C3) at (10,1) {\footnotesize \color{white} 4};
\node [node3] (C4) at (10,2.5) {\footnotesize \color{white} 3};

  \draw [-latex,draw=red, very thick](C2) -- (C1);
  \draw [-latex,draw=red, very thick](C3) -- (C2);
  \draw [-latex,draw=red, very thick](C4) -- (C3);
  \draw [-latex,draw=red, very thick](C1) -- (C4);

  \tikzstyle{node4} = [circle,inner sep=1pt,draw=blue,fill=blue,opacity=0.5];

\node [node4] (C1423) at (0.803278689*\mainfact,0.68852459*\mainfact,0.344262295*\mainfact) {};

\end{tikzpicture}

    \caption{The set of efficient weight vectors (the union of $\mathcal{T}_1, \mathcal{T}_2$ and $\mathcal{T}_3$) for the simple perturbed example, where the $(1,2,3)$ and $(1,2,4)$ triads, and the $(1,4,2,3,1)$ $4$-cycle are consistent. An interactive version is available at https://www.geogebra.org/m/psetpjfx.}
    \label{SimplePerturbed}
\end{figure}

\setcounter{appendbfigure}{6}
\end{example}

\setcounter{appendbexample}{5}
\begin{example} \label{consistent}
Figure~\ref{Consistent} shows the set of efficient weight vectors for the following consistent PCM:

\[ \mathbf{A''} =
\begin{pmatrix}
1 & \color{orange} 5/2 & 5 & 7 \\
\color{orange}2/5 & 1 & 2 & \color{blue} 14/5 \color{black} \\
1/5 & 1/2 & 1 & \color{green}7/5 \\
1/7 & \color{blue}5/14 \color{black} & \color{green}5/7 & 1
\end{pmatrix}.
\]

Here the original example is only modified by three elements (and their reciprocals), which are highlighted by orange, blue and green colors in the matrix. All the triads and the $4$-cycles are consistent, thus all three tetrahedra correspond to a single point in Figure~\ref{Consistent}.

\renewcommand{\thefigure}{B\arabic{appendbfigure}}
\begin{figure}[H]
    \centering

\begin{tikzpicture}[line join = round, line cap = round, scale=0.87]

\pgfmathsetmacro{\mainfact}{8};

\coordinate [label=above left:\text{$V_1=(w_1=1,w_2=0,w_3=0,w_4=0)$}] (A) at (1*\mainfact,1*\mainfact,0);
\coordinate [label=below:\text{$V_2=(w_1=0,w_2=1,w_3=0,w_4=0)$}] (B) at (1*\mainfact,0,1*\mainfact);
\coordinate [label=left:\text{$V_3=(w_1=0,w_2=0,$}] (C) at (0,1*\mainfact,1*\mainfact);
\coordinate [label=left:\text{$w_3=1,w_4=0)$}] (CC) at (0,0.93*\mainfact,1*\mainfact);
\coordinate [label=left:\text{$V_4=(w_1=0,w_2=0,w_3=0,w_4=1)$}] (D) at (0,0,0);

\draw[->] (0,0) -- (10,0,0) node[right] {$x$};
\draw[->] (0,0) -- (0,10,0) node[above] {$y$};
\draw[->] (0,0) -- (0,0,10) node[below left] {$z$};

\foreach \i in {A,B,C,D}
    \draw[dashed,very thick] (0,0)--(\i);
\draw[-,very thick, opacity=.5] (A)--(D)--(B)--cycle;
\draw[-,very thick, opacity=.5] (A) --(D)--(C)--cycle;
\draw[-,very thick, opacity=.5] (B)--(D)--(C)--cycle;
\draw[-,very thick] (A)--(B);
\draw[-,very thick] (A)--(C);
\draw[-,very thick] (A)--(D);

\coordinate (I) at (0.803278689*\mainfact,0.68852459*\mainfact,0.344262295*\mainfact);
\coordinate (J) at (0.803278689*\mainfact,0.68852459*\mainfact,0.344262295*\mainfact);
\coordinate (K) at (0.803278689*\mainfact,0.68852459*\mainfact,0.344262295*\mainfact);
\coordinate (L) at (0.803278689*\mainfact,0.68852459*\mainfact,0.344262295*\mainfact);

\draw[-,very thick, fill = blue!30, opacity =.5] (I)--(J)--(K)--cycle;
\draw[-,very thick, fill = blue!30, opacity =.5] (I) --(K)--(L)--cycle;
\draw[-,very thick, fill = blue!30, opacity =.5] (J)--(K)--(L)--cycle;

\coordinate (M) at (0.803278689*\mainfact,0.68852459*\mainfact,0.344262295*\mainfact);
\coordinate (N) at (0.803278689*\mainfact,0.68852459*\mainfact,0.344262295*\mainfact);
\coordinate (O) at (0.803278689*\mainfact,0.68852459*\mainfact,0.344262295*\mainfact);
\coordinate (P) at (0.803278689*\mainfact,0.68852459*\mainfact,0.344262295*\mainfact);

\draw[-,very thick, fill = red!30, opacity =.5] (M)--(N)--(O)--cycle;
\draw[-,very thick, fill = red!30, opacity =.5] (M) --(O)--(P)--cycle;
\draw[-,very thick, fill = red!30, opacity =.5] (N)--(O)--(P)--cycle;

\coordinate (E) at (0.803278689*\mainfact,0.68852459*\mainfact,0.344262295*\mainfact);
\coordinate (F) at (0.803278689*\mainfact,0.68852459*\mainfact,0.344262295*\mainfact);
\coordinate (G) at (0.803278689*\mainfact,0.68852459*\mainfact,0.344262295*\mainfact);
\coordinate (H) at (0.803278689*\mainfact,0.68852459*\mainfact,0.344262295*\mainfact);

\draw[-,very thick, fill=green!30, opacity=.5] (E)--(H)--(F)--cycle;
\draw[-,very thick, fill=green!30, opacity=.5] (E) --(H)--(G)--cycle;
\draw[-,very thick, fill=green!30, opacity=.5] (F)--(H)--(G)--cycle;

\tikzstyle{node1} = [circle,inner sep=2pt,draw=forestgreen,fill=forestgreen,thick];

\node [node1] (A1) at (8.5,7.5) {\footnotesize \color{white} 1};
\node [node1] (A2) at (10,7.5) {\footnotesize \color{white}  2};
\node [node1] (A3) at (10,6) {\footnotesize \color{white} 3};
\node [node1] (A4) at (8.5,6) {\footnotesize \color{white} 4};

  \draw [latex-latex,draw=forestgreen, very thick](A1) -- (A2);
  \draw [latex-latex,draw=forestgreen, very thick](A2) -- (A3);
  \draw [latex-latex,draw=forestgreen, very thick](A3) -- (A4);
  \draw [latex-latex,draw=forestgreen, very thick](A4) -- (A1);

  \tikzstyle{node2} = [circle,inner sep=2pt,draw=blue,fill= blue,thick];

\node [node2] (B1) at (8.5,5) {\footnotesize \color{white} 1};
\node [node2] (B2) at (8.5,3.5) {\footnotesize \color{white} 3};
\node [node2] (B3) at (10,3.5) {\footnotesize \color{white} 2};
\node [node2] (B4) at (10,5) {\footnotesize \color{white} 4};

  \draw [latex-latex,draw=blue, very thick](B2) -- (B1);
  \draw [latex-latex,draw=blue, very thick](B3) -- (B2);
  \draw [latex-latex,draw=blue, very thick](B4) -- (B3);
  \draw [latex-latex,draw=blue, very thick](B1) -- (B4);

\tikzstyle{node3} = [circle,inner sep=2pt,draw=red,fill= red,thick];

\node [node3] (C1) at (8.5,2.5) {\footnotesize \color{white} 1};
\node [node3] (C2) at (8.5,1) {\footnotesize \color{white} 2};
\node [node3] (C3) at (10,1) {\footnotesize \color{white} 4};
\node [node3] (C4) at (10,2.5) {\footnotesize \color{white} 3};

  \draw [latex-latex,draw=red, very thick](C2) -- (C1);
  \draw [latex-latex,draw=red, very thick](C3) -- (C2);
  \draw [latex-latex,draw=red, very thick](C4) -- (C3);
  \draw [latex-latex,draw=red, very thick](C1) -- (C4);

  \tikzstyle{node4} = [circle,inner sep=1pt,draw=blue,fill=blue,opacity=0.5];




\tikzstyle{node6} = [circle,inner sep=1pt,fill=green];

\node [node6] (C1234) at (0.803278689*\mainfact,0.68852459*\mainfact,0.344262295*\mainfact){};

\end{tikzpicture}

    \caption{The set of efficient weight vectors (the union of $\mathcal{T}_1, \mathcal{T}_2$ and $\mathcal{T}_3$) for the consistent example.}
    \label{Consistent}
\end{figure}

\setcounter{appendbfigure}{7}
\end{example}

\end{document}